\documentclass[a4paper,12pt]{article}

\usepackage{latexsym}
\usepackage{changebar} 
\usepackage{comment}   
\outerbarstrue         

\newcommand{\be}{\begin{equation}}      
\newcommand{\ee}{\end{equation}}        
\newcommand{\bern}{\begin{eqnarray*}}   
\newcommand{\eern}{\end{eqnarray*}}     
\newcommand{\beqp}{\begin{eqproof}}     
\newcommand{\eeqp}{\end{eqproof}}       

\newcommand{\bt}{\begin{teorema}}       
\newcommand{\et}{\end{teorema}}         
\newcommand{\bd}{\begin{definizione}}   
\newcommand{\ed}{\end{definizione}}     
\newcommand{\bc}{\begin{corollario}}    
\newcommand{\ec}{\end{corollario}}      
\newcommand{\bp}[1]{\noindent{\bf Proof #1.} } 
\newcommand{\ep}{\hfill$\Box$\par\medskip }   

\newcommand{\bl}{\begin{lemma}}         
\newcommand{\el}{\end{lemma}}           
\newcommand{\boss}{\begin{osservazione}\rm}     
\newcommand{\eoss}{\end{osservazione}}  
\newcommand{\refe}[1]{(\ref{#1})}       

\font\corsivo=rsfs10 at 12pt
\font\doppio=msbm10 at 12pt
\font\scdoppio=msbm8
\font\fs=cmr9

\newcommand{\C}{\mbox{\corsivo C}}      
\newcommand{\R}{\hbox{\doppio R}}       
\newcommand{\RR}{\hbox{\scdoppio R}}    
\newcommand{\nat}{\hbox{\doppio N}}     
\newcommand{\G}{\hbox{\doppio G}}     
\newcommand{\GG}{\hbox{\scdoppio G}}    

\newcommand{\norm}[1]{\left\|#1\right\|}                
\newcommand{\abs}[1]{\left|#1\right|}   
\newcommand{\decl}{:=}                  

\newcommand{\derivpar}[1]{\frac{\partial}{\partial #1}}

\newtheorem{teorema}{Theorem}[section]
\newtheorem{definizione}[teorema]{Definition}
\newtheorem{corollario}[teorema]{Corollary}
\newtheorem{lemma}[teorema]{Lemma}
\newtheorem{proposizione}[teorema]{Proposition}
\newtheorem{example}[teorema]{Example}

\newtheorem{osservazione}[teorema]{Remark}

\def\acknowledgment{\goodbreak\subsection*{Acknowledgment}
\bgroup \footnotesize}

\def\acknowledgments{\goodbreak\subsection*{Acknowledgments}
\bgroup \footnotesize}

\def\endacknowledgment{\vskip1sp\egroup}
\def\endacknowledgments{\vskip1sp\egroup}

\newcommand{\diver}{\mathrm{div}}       

\newcommand{\Cdue}{\mbox{\corsivo C}^{\,2}}
\newcommand{\Cuno}{\mbox{\corsivo C}^{\,1}}

\newcommand{\h}{\hbox{\doppio H}}     
\newcommand{\hei}{\hbox{\doppio H}}     
\newcommand{\hh}{\hbox{\scdoppio H}^n}  

\newcommand{\RN}{\R^N}

\newcommand{\grl}{\nabla_{\!\!L}}               
\newcommand{\diverl}{\diver_{\!L}}      

\newcommand{\grh}{\nabla_{\!\!H}}               
\newcommand{\lh}{\Delta_H}              

\newcommand{\grgr}{\nabla_{\!\!\gamma}}               
\newcommand{\lgr}{\Delta_\gamma}              
\newcommand{\normg}[1]{[\![#1]\!]}   

\newcommand{\norg}{\normg\xi}

\title{Hardy Type Inequalities Related to\\
Degenerate Elliptic Differential Operators}
\author{{Lorenzo D'Ambrosio }\\
{\fs Dipartimento di Matematica, via E. Orabona, 4 --
  I-70125, Bari}\thanks{Address for correspondence.
  e-mail: dambros@dm.uniba.it}}
\date{}
\begin{document}
\maketitle
{\abstract {We prove some Hardy type
inequalities related to quasilinear second order degenerate elliptic
differential operators $L_pu\decl-\nabla_L^*(\abs{\nabla_Lu}^{p-2}\nabla_Lu)$. If $\phi$ is a 
positive weight 
such that $-L_p\phi\ge 0$, then the Hardy type inequality 
$$ c\int_\Omega \frac{\abs u^p}{\phi ^p}\abs{\grl \phi}^p\ d\xi
  \le \int_\Omega\abs{\grl u}^p  d\xi\qquad \left(u\in\Cuno_0(\Omega)\right)$$
holds. We find an explicit value of the
constant involved, which, in most cases, results optimal. As
particular case we derive Hardy inequalities for subelliptic
operators on Carnot Groups.

}

}


\section{Introduction}\label{sec:intrhi}

An $N$-dimensional generalization of the classical Hardy
inequality is the following
\be c \int_\Omega \abs{u}^p w^{-p} dx
        \le \int_\Omega \abs{\nabla u}^p dx,\qquad
        u\in\Cuno_0(\Omega),\label{hardyeu}\ee
where $p>1$, $\Omega\subset\RN$ and the weight $w$ is, for instance,
$w\decl\abs x$ or $w(x)\decl dist(x,\partial\Omega)$
(see for instance \cite{bar-fil-ter04, bre-mar97, dav95}
and the references therein).

A lot of efforts have been made to give explicit values of the constant
$c$, and even more, to find its best value $c_{n,p}$ (see e.g.
\cite{bar-fil-ter04, bre-mar97, dav95, dav-hin98, gaz-gru-mit04, mar-miz-pin98, mat-sob97, mat-sob98}).

The preeminent rule of the Hardy inequality in the study of linear and
nonlinear
partial differential equations is well-known.
For instance, let us consider the  linear initial value
problem
\be\left\{\begin{array}{ll} u_t-\Delta u=\lambda\frac{u}{\abs{x}^2},&
        x\in \R^n,\quad n\ge3,\quad t\in]0,T[,\quad\lambda\in\R,\cr
        u(x,0)=u_0(x),&x\in\R^n,\quad u_0\in L^2(\R^n),\quad u_0>0.\cr\end{array}\right.\label{pb}\ee
The problem \refe{pb} has a solution if and only if $\lambda\le (\frac{n-2}{2})^2=c_{n,2}$
(see \cite{bar-gol84} for more details).
In the last years this result has been extended in several directions
see e.g.
\cite{bre-cab98, bre-vaz97, gar-per98, gol-zha01, mit-poh01-1, mit-poh01b, vaz-zua00}.

In the Heisenberg group setting,
Garofalo and Lanconelli in \cite{gar-lan90},
Niu, Zhang and Wang in \cite{niu-zha-wan01} and the author in \cite{dam04-1}
proved, among other results,
the following Hardy type inequality related to the sub-Laplacian  $\Delta_H$ 
on the Heisenberg group $\h^n$:
\be c\int_{\hh}\frac{u^2}{\rho^2}\psi_H^2d\xi\le\int_{\hh}\abs{\nabla_H u}^2d\xi,
     \quad u\in\Cuno_0(\h^n\setminus\{0\})   \label{hardy-gar-lan} \ee
where $\nabla_H$ denotes the vector field associated to the real part of the Kohn Laplacian
($\Delta_H=\nabla_H\cdot\nabla_H$), $\rho$ and $\psi_H$ are respectively
a suitable distance from the origin and a weight function
such that $0\le\psi_H\le1$.

 Recently, in \cite{gol-zha01}, it has been pointed out that
the analogue problem of \refe{pb} involving the sub-Laplacian $\Delta_H$,
namely
\[\left\{\begin{array}{ll} u_t-\Delta_H u=\lambda\psi_H^2\frac{u}{\rho^2}&
        {\rm on\,} \R^{2n+1}\times]0,T[,\quad\lambda\in\R,\cr
        u(\cdot,0)=u_0(\cdot)& {\rm on\,} \R^{2n+1},\quad u_0\in L^2(\R^{2n+1}),\quad u_0>0,\cr
        \end{array}\right.\]
has a positive solution if and only if $\lambda\le c_{b,H}$, where $c_{b,H}$ is
the best constant in \refe{hardy-gar-lan}.

Similar results have been established for equations involving
the Baouendi-Grushin type operators
$\Delta_\gamma:=\Delta_x+|x|^{2\gamma}\Delta_y=\grgr\cdot\grgr$
(see \cite{kom05}).

Recently, in \cite{dam-mit-poh04} Mitidieri, Pohozaev and the author
among other results, find some conditions
on the functions $u$ and $f$, that assure the positivity of the solutions
of the  partial differential inequalities $-Lu\ge f(\xi,u)$ on $\RN$.
Here $L$ is a quite general linear second order differential operator,
namely, $Lu\decl -\nabla_L^*\cdot\nabla_Lu$, where $\nabla_L$
is a general vector field.
This class of operators include all previous cited operators
as well as the sub-Laplacian on Carnot groups.

Having in mind some extensions of the above results
in the setting of second order linear degenerate (or singular)
partial differential operators, it appears that an important
step towards this programme is to establish some
fundamental inequalities of Hardy type.

In this paper we shall prove some Hardy type inequalities associated
to the quasilinear operators
$$L_pu\decl-\nabla_L^*(\abs{\nabla_Lu}^{p-2}\nabla_Lu)\qquad (p>1).$$
Our principal result can be roughly described as follows:
let $\phi:\Omega\to\R$ be any positive weight, 
 for any $u\in\Cuno_0(\Omega)$ we have
$$ c\int_\Omega \frac{\abs u^p}{\phi ^p}\abs{\grl \phi}^p\ d\xi
  \le \int_\Omega\abs{\grl u}^p d\xi,$$
provided $-L_p\phi\ge 0$.

For this goal we shall mainly use a technique developed in \cite{dam04-2, dam04-1, mit00}.
An interesting outcome of this approach is that, in several cases,
one can  easily obtain the best constant.
Furthermore, our main results represent a generalization of
some results contained in \cite{bar-fil-ter03, bar-fil-ter04}.
Indeed, in those papers the authors deal with a very special case:
The usual Euclidean case where $\phi$ is a particular power of
the Euclidean distance from a given surface.
Whereas, in our approach, $\grl$ can be any quite general vector field
and $\phi$ any positive weight:
the generality of this approach being an important strength.
It is, in fact, to remark that this unifying method allows,
specializing the choice of $\phi$, to obtain
almost all the fundamental Hardy inequalities known 
in Euclidean and subelliptic settings as well as to yield
new Hardy type inequalities.
Moreover, let us to stress that our only hypothesis 
$-L_p\phi\ge 0$ plays a relevant role in order to establish
that the best constant is not achieved.

We pay particular attention to the following special cases of $L_p$:
the Grushin type operators, the Heisenberg-Greiner operators
and the sub-Laplacian on Carnot groups (see section 3).
Specializing the function $\phi$, we get more concrete Hardy
type inequalities for these operators with explicit values of the constants
involved, which result the best possible in almost all the considered cases.

\section{Main results}\label{sec:prehi}
The aim of this section is to present some preliminary results and derive
some Hardy type inequalities related to a general vector field.

In this paper $\nabla$ stands for the usual gradient in $\RN$.
We indicate with $I_k$ and with $\abs\cdot$ respectively
the identity matrix of order $k$ and the Euclidean norm.

Let $\mu\decl(\mu_{ij})$, $i=1,\dots,l$, $j=1,\dots,N$
be a matrix with continuous entries $\mu_{ij}\in\C(\R^N)$.
Let $X_i$ ($i=1,\dots,l$) be defined as
\be X_i\decl\sum_{j=1}^N \mu_{ij}(\xi)\frac{\partial}{\partial \xi_j}
   \label{mu}\ee
and let $\grl$ be the vector field defined by
$$\grl\decl (X_1,\dots,X_l)^T=\mu\nabla.$$
Assuming that for $i=1,\dots,l$ and $j=1,\dots,N$ the derivative
$\frac{\partial}{\partial\xi_j}\mu_{ij}\in\C(\RN)$, 
we set
$$X_i^*\decl-\sum_{j=1}^N \frac{\partial}{\partial\xi_j}\mu_{ij}(\xi)\cdot$$
the formal adjoint of $X_i$ and $\grl^*\decl(X_1^*,\dots,X_l^*)^T$.

For any vector field $h=(h_1,\dots,h_l)^T\in\Cuno(\Omega,\R^l)$, we shall use
the following notation
\[\diverl(h)\decl\diver(\mu^T h),\]
that is
\[ \diverl(h)=-\sum_{i=1}^l X_i^*h_i=-\grl^*\cdot h.\]
In what follows $L$ stands for the linear second order differential operator defined by
$$L\decl \diverl(\grl)=-\sum_{i=1}^l X_i^*X_i=-\grl^*\cdot\grl$$
and for $p>1$, with $L_p$ we denote the quasilinear operator
$$L_p(u)\decl \diverl(\abs{\grl u}^{p-2}\grl u)=
   -\sum_{i=1}^l X_i^*(\abs{\grl u}^{p-2} X_iu)=
   -\grl^*\cdot(\abs{\grl u}^{p-2}\grl u).$$

\begin{example} Let $l<N$ be a positive natural number and let $\mu^l$
be the matrix defined as
  $$\mu^l\decl\pmatrix{ I_l & 0 }.$$
  The corresponding vector field $\nabla^l$ results to be the usual gradient
  acting only on the first $l$ variables
  $\nabla^l=(\frac{\partial}{\partial\xi_1},\frac{\partial}{\partial\xi_2},\dots,\frac{\partial}{\partial\xi_l})$.
  It is clear that $\nabla^N=\nabla$.
  The corresponding quasilinear operator $L_p$ is the usual
  $p-$Laplacian acting on the first $l$ variables of  $\RN$.
\end{example}

\begin{example}(Baouendi-Grushin type operator)
  Let $\RN$ be splitted in $\xi=(x,y)\in\R^n\times\R^k$.
  Let $\gamma\ge 0$ and let $\mu$ be the following matrix
  \be 
    \pmatrix{ I_n & 0 \cr 0 & \abs x^\gamma I_k }.\label{matgrushin}\ee
  The corresponding vector field is  $\grgr=(\nabla_x,\abs x^\gamma\nabla_y)$
  and the linear operator $L$ is the so-called Baouendi-Grushin operator
  $L=\Delta_x+\abs x^{2\gamma}\Delta_y$.

  Notice that if $k=0$ or $\gamma=0$, then $L$ and $L_p$ coincide  respectively
  with the usual Laplacian operator and $p$-Laplacian operator.
\end{example}

\begin{example}(Heisenberg gradient) Let $\xi=(x,y,t)\in\R^n\times\R^n\times\R=\h^n(=\RN)$ and
  let $\mu$ be defined as
  \[ 
    \left(\matrix{ I_n & 0 & 2y \cr 0 & I_n & -2x }\right).\]
  The corresponding vector field $\grh$ is the Heisenberg gradient on
  the Heisenberg group $\h^n$.

  This is the simplest case of a more general setting: the Carnot groups.
  More details are given in Section \ref{sec:carnot}.
\end{example}

 \begin{example}(Heisenberg-Greiner operator)
   Let $\xi=(x,y,t)\in\R^n\times\R^n\times\R$, $r\decl\abs{(x,y)}$,
   $\gamma\ge1$ and let $\mu$ be defined as
  \be 
    \left(\matrix{ I_n & 0 & 2\gamma yr^{2\gamma-2} \cr
      0 & I_n & -2\gamma xr^{2\gamma-2} }\right).\label{matgrein} \ee
  The corresponding vector fields are
  $X_i=\derivpar{x_i}+2\gamma y_ir^{2\gamma-2}\derivpar{t}$,
  $Y_i=\derivpar{y_i}-2\gamma x_ir^{2\gamma-2}\derivpar{t}$ for $i=1,\dots,n$.

  For $\gamma=1$ $L$ is the sub-Laplacian $\lh$ on
  the Heisenberg group $\h^n$.
  If $\gamma=2,3,\dots$, $L$ is a Greiner operator
  (see \cite{gre79}).
\end{example}

Let $A$ be an open subset of $\R^N$ with Lipschitz boundary $\partial A$ and
let $\hat h\in\Cuno(\overline A,\R^l)$ be
a vector field. By the divergence theorem we have
\[ \int_A \diverl{\hat h}d\xi =
\int_A \diver{(\mu^T \hat h)}d\xi =
        \int_{\partial A} \hat h\cdot\mu \nu d\Sigma =
        \int_{\partial A} \hat h\cdot \nu_L d\Sigma,\]
where $\nu_L\decl\mu\nu$, and $\nu$ denotes the exterior normal
at point $\xi\in\partial A$.
If $\hat h$ has the form $\hat h= fh$ with
$f\in \Cuno(\overline A)$ and $h\in\Cuno(\overline A,\R^l)$, then
\be\int_A f\diverl{h}d\xi +
        \int_{A}\grl f\cdot h d\xi =
        \int_{\partial A} f  h\cdot \nu_L d\Sigma\label{for:GGH}.\ee
Moreover, if $h=\grl u$ with $u\in \Cdue(\overline A)$, then
(\ref{for:GGH}) yields the Gauss--Green formula
\[\int_A f Lu d\xi + \int_A \grl f\cdot\grl u d\xi =
        \int_{\partial A}f\grl u\cdot \nu_L d\Sigma.\]

Let $g\in\Cuno(\R)$ be such that $g(0)=0$ and let $\Omega\subset\R^N$ be open.
For every vector field $h\in\Cuno(\overline A,\R^l)$ and any compactly
supported function $u\in\Cuno_0(\Omega)$,
choosing $f\decl g(u)$ in (\ref{for:GGH}), we obtain
\be \int_\Omega g(u)\diverl{h}d\xi = - \int_{\Omega} g'(u)\grl u\cdot h d\xi.\label{id1} \ee
Let  $h\in L^1_{loc}(\Omega,\R^{l})$ be a vector field.
As usual, we define the distribution $\diverl h$
using the formula (\ref{id1}) with $g(s)=s$.
If in (\ref{id1}) we chose $g(t)=\abs{t}^p$ with $p>1$,
then for every $u\in\Cuno_0(\Omega)$  we have
\be\int_\Omega \abs{u}^p \diverl{h}d\xi = - p \int_{\Omega} \abs{u}^{p-2}u
        \grl u\cdot h d\xi.\label{id2}\ee

Let $h\in L^1_{loc}(\Omega,\R^{l})$ be a  vector field and
let $A\in L^1_{loc}(\Omega)$ be a function.
In what follows we write $A\le \diverl h$ meaning that the
inequality holds in distributional sense, that is for every
 $\phi\in\Cuno_0(\Omega)$ such that $\phi\ge 0$, we have
$$\int_\Omega \phi A\, d\xi \le \int_\Omega \phi \diverl{h}d\xi =
  - \int_{\Omega} \grl \phi\cdot h\,d\xi.$$

Identities \refe{id1} and \refe{id2} play an important role in the proof of
the following Hardy type inequalities and the Poincar\`e inequality too.

\bt\label{teo:hardyh} Let $p>1$.
  Let $h\in L^1_{loc}(\Omega,\R^{l})$ be a vector field and let
  $A_h\in L^1_{loc}(\Omega)$ be a nonnegative function such that
  $A_h\le \diverl h$
  and $\abs{h}^pA_h^{1-p} \in L^1_{loc}(\Omega)$.
  Then for every $u\in\Cuno_0(\Omega)$, we have
\be \int_\Omega \abs{u}^p A_h\ d\xi \le  p^p \int_{\Omega}
   \frac{\abs{h}^p}{A_h^{(p-1)}}\abs{\grl u}^{p} d\xi. \label{dis:hardyh}\ee
\et

\bp{} We note that the right hand side of \refe{dis:hardyh} is finite since
        $u\in \Cuno_0(\Omega)$.
Using the identity (\ref{id2}) and H\"older inequality we obtain
\bern \int_\Omega \abs{u}^p A_h d\xi &\le& \int_\Omega \abs{u}^p \diverl{h}d\xi
        \le p \int_{\Omega} \abs{u}^{p-1} \abs{ h} \abs{\grl u} d\xi\\
        &=& p \int_{\Omega} \abs{u}^{p-1}A_h^{(p-1)/p}
        \frac{\abs{h}}{A_h^{(p-1)/p}} \abs{\grl u} d\xi\\
        &\le& p\left(\int_{\Omega} \abs{u}^{p}A_hd\xi\right)^{(p-1)/p}
        \left(\int_{\Omega}\frac{\abs{ h}^p}{A_h^{p-1}} \abs{\grl u}^p d\xi\right)^{1/p}.\eern
This completes the proof.
\ep

Specializing the vector field $h$ and the function $A_h$,
we shall deduce  from \refe{dis:hardyh}
some concrete inequalities of Hardy type.

\boss
Setting $A_h=\diverl h$ in (\ref{dis:hardyh}), we have
\be \int_\Omega \abs{u}^p \diverl h\ d\xi \le  p^p \int_{\Omega}
   \frac{\abs{h}^p}{|\diverl h|^{(p-1)}}\abs{\grl u}^{p} d\xi. \label{dis:predav}\ee

Acting as Davies and Hinz in \cite{dav-hin98}, the choice
$h:=\grl V$ with $V$ such that  $LV>0$, yields
\be \int_\Omega \abs{u}^p \abs{LV} d\xi \le  p^p \int_{\Omega}
   \frac{\abs{\grl V}^p}{\abs{L V}^{(p-1)}}\abs{\grl u}^{p} d\xi. \label{dis:dav}\ee
In order to state a Hardy inequality,
now the problem is to find a suitable function $V$.
In the Euclidean setting for $1<p<N$, choosing $V(\xi)=\abs{\xi}^{2-p}$
if $1<p<2$, $V(\xi)=\ln\abs\xi$ if $p=2$ and $V(\xi)=-\abs{\xi}^{2-p}$
if $2<p<N$,
we obtain the Hardy inequality (\ref{hardyeu}) with
$w(\xi)=\abs{\xi}$.

Another strategy is to chose the vector field $h$ as
$h=\abs{\grl V}^{p-2}\grl V$ with $V$ such that  $L_pV>0$. Thus, we have
\be \int_\Omega \abs{u}^p \abs{L_pV} d\xi \le  p^p \int_{\Omega}
   \frac{\abs{\grl V}^{p(p-1)}}{\abs{L_p V}^{(p-1)}}\abs{\grl u}^{p} d\xi. \label{dis:dav+}\ee
Hence, in the Euclidean setting  for $1<p<N$, choosing
$V(\xi)=\ln\abs \xi$ we reobtain the inequality (\ref{hardyeu}) with
$w(\xi)=\abs{\xi}$.

In order to obtain the classical Hardy inequalities in Euclidean setting,
these strategies are equivalent.
This equivalence is basically due to the fact that$\abs{\nabla\abs \xi}=1$
for $\xi\neq 0$.
The latter approach is slightly more simple:
the choice of $V$ is independent of $p$.
Moreover, it turned out to be more fruitful in the Heisenberg group
and in the Grushin plane settings (see \cite{dam04-1, dam04-2})
as well as in our more general framework.
\eoss

Let $d:\Omega\to\R$ be a nonnegative non constant measurable function.
In order to state Hardy inequalities involving the weight $d$,
the basic assumption we made on $d$ is that, for $\alpha\neq 0$,
$d^\alpha$ is a one side weak solution of $-L_p(u)=0$, that is
$d^\alpha$ is super-$L_p$-harmonic or sub-$L_p$-harmonic in weak sense.
Namely, let $\alpha,\beta\in\R$, $\alpha\neq 0$, requiring
\be d^{(\alpha-1)(p-1)}\abs{\grl d}^{p-1}\in L^1_{loc}(\Omega),\label{h1}\ee
we assume that
\be -L_p(d^\alpha)\ge 0\quad{\rm on\ }\Omega\qquad[\mathrm{resp.} \le0]\label{eq:har}\ee
in weak sense, that is for every nonnegative $\phi\in\Cuno_0(\Omega)$
we have
\be \int_\Omega \abs{\grl d^\alpha}^{p-2} \grl d^\alpha\cdot\grl \phi =
  \alpha\abs{\alpha}^{p-2}\int_\Omega d^{(\alpha-1)(p-1)}\abs{\grl d}^{p-2}\grl d\cdot\grl \phi\ge 0\
  [\mathrm{resp.} \le 0] \label{dis:sup}\ee
and
\be \alpha[(\alpha-1)(p-1)-\beta-1]>0,\qquad[\mathrm{resp.} <0].\label{h3} \ee

Gluing together the above conditions, we assume that
\be -L_p(cd^\alpha)\ge 0\quad{\rm on\ }\Omega\label{eq:shar}\ee
in weak sense, where $c\decl \alpha[(\alpha-1)(p-1)-\beta-1] $

\bt\label{th:hm} Assume that (\ref{h1}) and (\ref{eq:shar}) hold.
  Let $\beta \in \R$ be such that
  \begin{eqnarray}
    d^{\beta}\abs{\grl d}^{p}\in L^1_{loc}(\Omega),\label{h4}\\
    d^{\beta+p}\in L^1_{loc}(\Omega).\label{h5}
  \end{eqnarray}
  For every function $u\in\Cuno_0(\Omega)$, we have
  \be (c_{\alpha,\beta,p})^p
    \int_\Omega {\abs u^p}{d^{\beta}} \abs{\grl d}^{p} d\xi
    \le\int_\Omega d^{\beta+p}\abs{\grl u}^p d\xi,\label{dis:hargen}\ee
  where $c_{\alpha,\beta,p}\decl{\abs{(\alpha-1)(p-1)-\beta-1}}/{p}$.

  In particular, if $-L_p(d^\alpha)\ge 0$, then for every function $u\in\Cuno_0(\Omega)$ we have
  \be \left(\frac{\abs{\alpha}(p-1)}{p}\right)^p\int_\Omega
    \frac{\abs u^p}{d^{p}}\abs{\grl d}^{p} d\xi
      \le\int_\Omega \abs{\grl u}^p d\xi\label{dis:hargenpart}\ee
 provided $d^{-p}\abs{\grl d}^{p}\in L^1_{loc}(\Omega)$.
\et

\boss\label{rm:sharp} In most examples we shall deal with, the
  constant $c_{\alpha,\beta,p}^p$, yielded by applying Theorem
    \ref{th:hm}, results to be sharp.
    We shall now indicate an argument that can be used to prove the sharpness
    of the constant $c_{\alpha,\beta,p}^p$ involved in the inequality
    of Theorem \ref{th:hm}. Let $c_b(\Omega)$ be the best constant
    in (\ref{dis:hargen}). It is clear that $c_b(\Omega)\ge
    c_{\alpha,\beta,p}^p$.
    We shall assume that the hypotheses of
    Theorem \ref{th:hm} are satisfied and that
    there exists $s>0$ such that $\Omega^s\decl d^{-1}(]-\infty,s[)$
    and $\Omega_s\decl d^{-1}(]s,+\infty[)$ are not empty open subsets
    of $\Omega$ with piecewise regular boundaries.

  We assume that there exists $\epsilon_0>0$ such that for every
  $\epsilon\in ]0,\epsilon_0[$ there hold
\be 0<\int_{d<s} d^{c(\epsilon)p+\beta} \abs{\grl d}^p<+\infty,
  \ \ 0<\int_{d>s} d^{-c(\epsilon)p+\beta} \abs{\grl d}^p<+\infty,\label{int:sharp}
\ee
  where
\be\label{eq:cepsilon}
  c(\epsilon)\decl \frac{\abs{(\alpha-1)(p-1)-\beta-1}+\epsilon}{p}=
                c_{\alpha,\beta,p}+\frac\epsilon p .\ee

By rescaling argument, we can assume that $s=1$.
   Let $\epsilon\in ]0,\epsilon_0[$ and let $v:\Omega\to\R$ be defined as
   \be\label{eq:sharp} v(\xi)\decl \left\{\begin{array}{ll}
       d^{c(\epsilon)}(\xi) &\mathrm{if}\ d(\xi)\le 1,\\
       d^{-c(\epsilon/2)}(\xi)  &\mathrm{if}\ d(\xi) >1.
     \end{array}\right.\ee
   By hypothesis, $\int_\Omega v^p d^\beta\abs{\grl d}^p$ is finite.
   Thus, we have
   \begin{eqnarray*}
     c(\epsilon)^p \int_\Omega v^p d^\beta\abs{\grl d}^p&=&
     c(\epsilon)^p \int_{d<1} d^{\beta+p} d^{(c(\epsilon)-1)p}\abs{\grl d}^p
+c(\epsilon)^p \int_{d>1} d^{\beta+p} d^{(-c(\epsilon/2)-1)p}
         \abs{\grl d}^p\\
     &=&\int_{d<1} d^{\beta+p} \abs{\grl v}^p+
     (\frac{c(\epsilon)}{c(\epsilon/2)})^p\int_{d>1}d^{\beta+p}\abs{\grl v}^p\\
     &=&\int_{\Omega} d^{\beta+p} \abs{\grl v}^p+
     (\frac{c(\epsilon)^p}{c(\epsilon/2)^p}-1) \int_{d>1} d^{\beta+p} \abs{\grl v}^p.
   \end{eqnarray*}
   Observing that  $c(\epsilon)>c(\epsilon/2)$, we get
   \be\label{dis:over} c(\epsilon)^p \int_\Omega v^p d^\beta\abs{\grl d}^p>\int_{\Omega} d^{\beta+p} \abs{\grl v}^p,\ee
  the converse of the Hardy inequality.

  Now we assume that the Hardy inequality (\ref{dis:hargen}) holds for the
  function $v$ defined in (\ref{eq:sharp}).
  From (\ref{dis:over}) we deduce $c(\epsilon)^p>c_b(\Omega)$.
  Letting $\epsilon\to 0$, we get $c_{\alpha,\beta,p}^p\ge c_b(\Omega)$
  and hence the claim.
\eoss

The question of the existence of functions that realize the best constant
arises. In such a general framework 
a unique answer cannot be given. Indeed, even in the Euclidean setting several cases
occur. Let $p=2$, let $d_1(\cdot)\decl\abs\cdot$ be the Euclidean distance from the origin,
and let $d_2(\cdot)\decl dist(\cdot,\partial\Omega)$ be the distance from the boundary 
of a given domain $\Omega$. 
If $\Omega\subset\RN$ ($N\ge 3$) is a ball centered at the origin, then the best constants
in the Hardy inequality (\ref{dis:hargenpart}) related to $d_1$ and $d_2$ are not
achieved. On the other hand, there exist smooth bounded domains $\Omega$ such that
the best constant in the inequality related to $d_1$ is not achieved and 
the best constant in the inequality related
to $d_2$ is achieved (see \cite{ mar-miz-pin98, mat-sob97}). 
Anyway, some steps in this direction can be done even in our general framework.
For the sake of simplicity, we shall focus our attention on the inequality (\ref{dis:hargenpart}).

Therefore, under the same hypotheses of Theorem \ref{th:hm} we assume that
$ -L_p(d^\alpha)\ge 0$ on $\Omega$ in weak sense, that
$(\int_\Omega\abs{\grl u}^p d\xi)^{1/p}$ is a norm and that
$D^{1,p}_L(\Omega)$, the closure of $\C_0^{\,\infty}(\Omega)$ in that
norm, is well defined.
We denote by
 $c_b(\Omega)$ the best constant in (\ref{dis:hargenpart}), namely
\be c_b(\Omega)\decl\inf_{u\in D^{1,p}_L,u\neq0}\frac{\int_\Omega \abs{\grl u}^p d\xi}
{\int_\Omega {\abs u^p}{d^{-p}}\abs{\grl d}^{p} d\xi}
\label{def:bc}.
\ee
\bt\label{th:achieve} Under the above hypotheses we have:
\begin{enumerate}
\item If $d^{\alpha \frac{p-1}{p}}\in D^{1,p}_L(\Omega)$, then 
    $c_b(\Omega)=({\abs\alpha \frac{p-1}{p}})^p$ and $d^{\alpha \frac{p-1}{p}}$ is a minimizer.
\item  If $d^{\alpha \frac{p-1}{p}}\not\in D^{1,p}_L(\Omega)$, $p\ge2$, $\abs{\grl d}\neq0$ a.e.
    and $c_b(\Omega)=({\abs\alpha \frac{p-1}{p}})^p$ then the best constant 
	$c_b(\Omega)$ is not achieved.
\end{enumerate}
\et

\boss In all the examples we shall deal with in the last section,
it is possible to apply Theorem \ref{th:achieve} and, hence,
for $p\ge2$ the best constants mentioned in all the theorems of 
Section \ref{sec:apply} are not achieved.
\eoss

\boss Let us to consider the special case of $\grl=\nabla$, the usual
  Euclidean gradient, $d$ is the Euclidean distance from a given regular
  surface $K$ of codimension $k$ ($1\le k\le N$), $\alpha=\frac{p-k}{p-1}$
  and $\beta=-p$.
  In this case, replacing $\Omega$ with $\Omega\setminus K$, Theorem
  \ref{th:hm} assures that the inequality
   \be \left(\frac{\abs{p-k}}{p}\right)^p\int_\Omega
    \frac{\abs u^p}{d^{p}} d\xi
    \le\int_\Omega \abs{\nabla u}^p d\xi\label{dis:heucdist}\ee
  holds for every $u\in\Cuno_0(\Omega\setminus K)$ provided
  $-\Delta_p(d^\alpha)\ge 0$ on $\Omega\setminus K$.

  This particular case of Theorem \ref{th:hm} is contained in
  \cite{bar-fil-ter03, bar-fil-ter04}, where the authors also study the
  remainder terms for inequality (\ref{dis:heucdist}).

  The reader interested in the study of Hardy inequalities with remainder terms
  can refer to
  \cite{bar-fil-ter03, bar-fil-ter04, bre-mar97, bre-mar-sha00, gaz-gru-mit04}
  and the references therein for the Euclidean case and to
  \cite{dam04-1} for the case $\grl=\grh$, the Heisenberg gradient on the
  Heisenberg group.
\eoss

\bp{of Theorem \ref{th:hm}} We prove the thesis in the case
  $-L_p(d^\alpha)\ge 0$ and $c\decl \alpha[(\alpha-1)(p-1)-\beta-1]>0$.
  The alternative case is similar.

  Let $\varphi\in\Cuno_0(\Omega)$ be a nonnegative
  function.
  Choosing in (\ref{dis:sup}) $\phi:=d^{\beta+1-(\alpha-1)(p-1)}\varphi$,
  we have
  \be 0\le \alpha\int_\Omega d^{\beta+1}  \abs{\grl d}^{p-2}
       \grl d\cdot\grl \varphi -
     \alpha[(\alpha-1)(p-1)-\beta-1]\int_\Omega d^\beta \abs{\grl d}^{p}\varphi.\label{dis:ch} \ee
  Using H\"older inequality and hypotheses (\ref{h4}) and (\ref{h5}),
  it is immediate to check that the above integrals are finite.

  Let $h$ be the vector field defined by
  $h\decl -\alpha d^{\beta+1}\abs{\grl d}^{p-2} \grl d$ and
  let $A_h$ be the function defined as
  $A_h\decl\alpha[(\alpha-1)(p-1)-\beta-1] d^\beta \abs{\grl d}^{p}$.
  Thus, from (\ref{dis:ch}) and the fact that $c>0$, we obtain $\diverl h\ge A_h\ge 0$.
  Now we are in the position to apply Theorem \ref{teo:hardyh} and
  this concludes the proof.
\ep

\bp{of Theorem \ref{th:achieve}} 1) From (\ref{dis:hargenpart}), we have 
$c_b(\Omega)\ge({\abs\alpha \frac{p-1}{p}})^p$. It is immediate to check
that $u\decl d^{\alpha \frac{p-1}{p}}$ realizes the infimum in (\ref{def:bc}).

2) Let $u\in\C^{\,\infty}_0(\Omega)$. We define the functional $I$ as 
$$I(u)\decl \int_\Omega \abs{\grl u}^p d\xi-
\left(\frac{\abs{\alpha}(p-1)}{p}\right)^p\int_\Omega
    \frac{\abs u^p}{d^{p}}\abs{\grl d}^{p} d\xi.$$
The functional $I$ is non negative, and the best constant will be 
achieved, if and only if, 
$I(u)=0$ for some $u\in D^{1,p}_L(\Omega)$.

Let $v$ be the new variable $v\decl d^{-\gamma}u$ with $\gamma\decl\alpha \frac{p-1}{p}$.
By computation we have
\be \abs{\grl u}^2=\abs{\gamma}^2v^2d^{2\gamma-2}\abs{\grl d}^2+d^{2\gamma}\abs{\grl v}^2
 +2\gamma v d ^{2\gamma-1}(\grl d\cdot\grl v).\label{eq:gr2}
\ee
(If $d$ is not smooth enough, by standard argument one can consider $d_\epsilon$
a regularization of $d$ and after the computation taking the limit as 
$\epsilon \rightarrow 0$).

We remind that the inequality 
\be (\xi-\eta)^s\ge\xi^s-s\eta\xi^{s-1} \label{dis:ggm}\ee
holds for every $\xi,\eta,s\in\R$ with $\xi>0, \xi>\eta$ and $s\ge1$ (see
\cite{gaz-gru-mit04}). Applying (\ref{dis:ggm}) and (\ref{eq:gr2}) with
$s=p/2$, $\xi=\abs{\gamma}^2v^2d^{2\gamma-2}\abs{\grl d}^2$ and 
$\eta=-2\gamma v d ^{2\gamma-1}(\grl d\cdot\grl v)-d^{2\gamma}\abs{\grl v}^2$, we have
\bern \abs{\grl u}^p&\ge&\abs{\gamma}^p v^pd^{\gamma p} d^{-p}\abs{\grl d}^p + p\abs\gamma^{p-2}\gamma\abs v^{p-2}vd^{(\alpha-1)(p-1)}\abs{\grl d}^{p-2}(\grl d\cdot\grl v)\\
&&\qquad+\frac p 2 \abs\gamma^{p-2}\abs v^{p-2} d^{(\alpha-1)(p-1)+1}\abs{\grl d}^{p-2}\abs{\grl v}^2.
\eern
Taking into account that $u\decl d^{\gamma}v$ we have
$$I(u)\ge I_1(v)+I_2(v)$$
where 
\bern I_1(v)&\decl& \int_\Omega p\abs\gamma^{p-2}\gamma\abs v^{p-2}vd^{(\alpha-1)(p-1)}\abs{\grl d}^{p-2}(\grl d\cdot\grl v)\ d\xi,\\
I_2(v)&\decl&\frac p 2 \abs\gamma^{p-2}\int_\Omega\abs v^{p-2} d^{(\alpha-1)(p-1)+1}\abs{\grl d}^{p-2}\abs{\grl v}^2 d\xi.
\eern
Re-arranging the expression in $I_1$ and  integrating by parts we obtain
\bern \lefteqn{I_1(v)=\left(\frac{p-1}{p}\right)^{p-1}\int_\Omega \left(\grl\abs v^p\cdot
   \abs{\grl d^\alpha}^{p-2}\grl d^\alpha\right) d\xi}\\
   &&=\left(\frac{p-1}{p}\right)^{p-1}\!\!\! \int_{\partial\Omega} \abs v^p
   \abs{\grl d^\alpha}^{p-2}\left(\grl d^\alpha\cdot\nu_L \right) d\Sigma+
   \left(\frac{p-1}{p}\right)^{p-1}\!\!\!\int_\Omega \abs v^p(-L_p(d^\alpha)) d\xi\ge0,
\eern
where we have used the fact that $v\in\C^{\,\infty}_0(\Omega)$ and the hypothesis
$-L_p(d^\alpha)\ge0$.
On the other hand we can rewrite $I_2$ as 
$$ I_2(v)=\frac2 p\abs\gamma^{p-2}\int_\Omega d^{(\alpha-1)(p-1)+1}\abs{\grl d}^{p-2}
\abs{\grl \abs v^{\frac p 2}}^2\ d\xi. $$
Thus, we conclude that for any $u\in D^{1,p}_L(\Omega)$
$$I(u)\ge\frac2 p\abs\gamma^{p-2}\int_\Omega d^{(\alpha-1)(p-1)+1}\abs{\grl d}^{p-2} 
\abs{\grl \abs v^{\frac p 2}}^2\ d\xi, $$
and this inequality implies the non existence of minimizers in $D^{1,p}_L(\Omega)$.
\ep

\medskip

Specializing the function $d$, we shall deduce  from
Theorem  \refe{th:hm} some concrete inequalities of Hardy type.
A first example is the following.
We assume that
there exists $m\in\nat$, $1\le m\le l$ such that  the matrix $\mu$
in (\ref{mu}) has the following form
 \be \mu\decl\pmatrix{ I_m & \mu_1 \cr 0  &  \mu_2 }\label{matspec}\ee
where  $\mu_1$ and  $\mu_2$ denote matrixes with $m\times(N-m)$ and
$(l-m)\times(N-m)$ continuous entries respectively
and $I_m$  stands for the identity
matrix of order $m$. Notice that this case occurs in all the
examples cited above.

Set $\eta\decl(\xi_1,\dots,\xi_m)$, $\tau\decl(\xi_{m+1},\dots,\xi_N)$
and let $v_p$ be defined for
$\eta\in\R^m\setminus\{0\}$ as
\begin{equation}
  \label{eq:grufun2}
  v_p(\eta)\decl\left\{
    \begin{array}{ll}
      \abs \eta^\frac{p-m}{p-1}& \mathrm{if}\ p\neq m,\\
      -\ln\abs \eta &  \mathrm{if}\ p= m.
    \end{array} \right.
\end{equation}
The function $v_p$ is $p$-harmonic on $\R^m\setminus\{0\}\times\R^{N-m}$
for the Euclidean $p$-Laplacian acting on the
$\eta$ variable $\Delta_{p,\eta}$ and hence also for the quasilinear operator $L_p$.
Moreover, there exists a constant $l_p\neq0$ such that
$$ -\Delta_{p,\eta}v_p=l_p\delta_0 \quad \mathrm{on}\ \R^m$$
in weak sense, where $\delta_0$ is the Dirac distribution at $0\in\R^m$
and $l_p>0$ if and only if $1<p\le m$.
These relations allow us to apply Theorem \ref{th:hm}.

\bt\label{teo:harspec} Assume that $\mu$ has the form (\ref{matspec})
  and let $\beta\in\R$ be fixed.
\begin{enumerate}
\item  Let $1<p<\infty $ and let $\Omega\subset\RN$ be an open set.
  If $m+\beta<0$, then we also require that  $\Omega\subset(\R^m\setminus\{0\})\times\R^{N-m}$.
  Then for every $u\in\Cuno_0(\Omega)$, we have
  \begin{equation}\label{eq:harspec}
   b_\beta^p\int_\Omega  {\abs{u(\eta,\tau)}^p}{\abs \eta^{ \beta}} d\eta d\tau
                 \le\int_{\Omega} \abs{\grl u(\eta,\tau)}^{p}\abs \eta^{p+\beta} d\eta d\tau,
  \end{equation}
  where $b_{\beta}\decl\frac{\abs{m+\beta}}{p}$.

  In particular, for every $u\in\Cuno_0(\Omega)$, we obtain
  \begin{equation}
    \left(\frac{\abs{m-p}}{p}\right)^p\int_\Omega\frac{\abs{u(\xi)}^p}{\abs\xi^p}d\xi
    \le\left(\frac{\abs{m-p}}{p}\right)^p\int_\Omega\frac{\abs{u(\eta,\tau)}^p}{\abs \eta^p}d\eta d\tau
       \le \int_{\Omega} \abs{\grl u(\xi)}^{p} d\xi. \label{dis:hardyspec2}
  \end{equation}
\item  Let $p = m>1$.
  Let $R>0$ and set $\Omega\decl\{\xi=(\eta,\tau)\in\R^m\times \R^{N-m},\ \abs \eta<R\}$.
  If $\beta<-1$, then for every $u\in\Cuno_0(\Omega)$, we have
 \begin{equation}\label{eq:hardyspecln}
   \tilde b_\beta^p\int_\Omega  \frac{\abs{u(\eta,\tau)}^p}{\abs \eta^p}
       (\ln\frac R{\abs \eta})^\beta d\eta d\tau
     \le\int_{\Omega} \abs{\grl u(\eta,\tau)}^{p} (\ln\frac R{\abs \eta})^{p+\beta} d\eta d\tau,
  \end{equation}
  where $\tilde b_\beta\decl\frac{\abs{\beta+1}}{p}$.

  In particular, for every $u\in\Cuno_0(\Omega)$, we obtain
 \begin{equation}\left(\frac{p-1}{p}\right)^p
   \int_\Omega  \frac{\abs{u(\eta,\tau)}^p}{(\abs \eta\ln(R/\abs \eta))^p}d\eta d\tau
     \le\int_{\Omega} \abs{\grl u(\xi)}^{p} d\xi.
  \end{equation}
\end{enumerate}
\et
\boss It is easy to check that the inequality
(\ref{eq:hardyspecln}) holds also for $\beta>-1$ provided the set
$\Omega$ is replaced by $\Omega\decl\{(\eta,\tau)\in\R^m\times
\R^{N-m},\ 0<\abs \eta<R\}$. \eoss
\bp{} Let $1<p<m$. We claim
   that the function $v_p$ is
   super-$L_p$-harmonic on $\RN$. Indeed, let $\phi\in\Cuno_0(\RN)$
   be a non negative function. Observing that
        $\abs{\grl v_p}=\abs{\nabla_\eta v_p}$, we have
   \be \int_{\RR^N} -L_p v_p \,\phi d\xi=
        \int_{\RR^{N-m}}d\tau\int_{\RR^m}d\eta (-\Delta_{p,\eta}v_p)\phi
        =l_p\int_{\RR^{N-m}}\phi(0,\tau)d\tau\ge 0
   \ee
   Analogously, one can prove that $v_p$ is super-$L_p$-harmonic
   when $p=m$ and sub-$L_p$-harmonic when $p>m$.

   First we consider the case $p\neq m$. We choose $d^\alpha=v_p$
   with $d(\xi)=d(\eta,\tau)=\abs \eta$ and $\alpha=\frac{p-m}{p-1}$.
   Observing that $\abs{\grl\abs\eta}=\abs{\nabla_\eta\abs\eta}=1$ a.e.
   and that the integrability conditions (\ref{h1}),(\ref{h4}),
   (\ref{h5}) are satisfied, applying Theorem \ref{th:hm}
   we get (\ref{eq:harspec}).

   Let $p=m>1$. The choices $d(\eta,\tau)=\ln\frac{R}{\abs\eta}$ and
   $\alpha=1$ in Theorem \ref{th:hm} yield the inequality
   (\ref{eq:hardyspecln}).

  Finally, we prove the missing inequality (\ref{eq:harspec}) when $p=m$.
  We consider the case $m+\beta>0$.
  The case $m+\beta<0$ is analogous and the case
  $m+\beta=0$ is trivial.
  Let $\sigma>0$ be such that $m+\beta-\sigma>0$.
  We chose  $d(\xi)=\abs\eta$, $\alpha=\frac{\sigma}{m-1}$.
  In this case it easy to check that $d^\alpha$ is sub-$L_m$-harmonic
  on $\RN$, that is
  $$-L_m(d^\alpha)=-\diverl \left((\frac{\sigma}{m-1})^{m-1}
    \abs\eta^{\sigma-m+1}
      \nabla_\eta \abs \eta \right)=-\frac{\sigma^m}{{(m-1)}^{m-1}}
    \frac{1 }{\abs\eta^{m-\sigma}}\le 0.$$
  The constant $c$ in (\ref{eq:shar}) is
  $c=\frac{\sigma}{m-1}\ \frac{\sigma-m-\beta}{m}<0$.
  Hence, we are in the position to apply Theorem \ref{th:hm}; thus, we derive
  the inequality
    $$ (\frac{{m+\beta-\sigma}}{p})^p \int_\Omega  {\abs{u}^p}{\abs \eta^{ \beta}} d\xi
                 \le\int_{\Omega} \abs{\grl u}^{p}\abs \eta^{p+\beta}
                 d\xi.$$
  Letting $\sigma\to 0$, we get the claim.
\ep
\boss In the case $\mu=I_N$, the vector field $\grl$ is the
  usual gradient $\nabla$. For $m<N$,  inequalities  of type
  (\ref{eq:harspec}) are already present in \cite{maz85} and in
  \cite{mit00}. Secchi, Smets and Willem in \cite{sec-sme-wil03} prove that
  the constant $b_\beta^p$
  is optimal when $m+\beta>0$ and $\Omega=\RN$ (see next section for further
  generalization in this direction).
\eoss


An immediate consequence of Theorem \ref{teo:harspec} is a
Poincar\`e inequality for the vector field $\grl$.
The claim easily follows from inequality (\ref{eq:harspec})
with $\beta=0$.

\bt\label{poin_l} Let $\Omega$ be an open subset of $\R^N$ bounded
in $\xi_1$ direction, that is, there exists $M>0$ such that for
every $\xi\in\Omega$ it results $|\xi_1|\le M$.
Assume that
the matrix $\mu$ has the form (\ref{matspec}).

\noindent Then, for every $u\in \Cuno_0(\Omega)$, we have
\[  c^p\int_\Omega \abs{u}^p d\xi
        \le\int_\Omega\abs{\grl u}^pd\xi,\]
with $c\decl\frac{1}{pM}$.
\et

In \cite{dam-mit-poh04} the authors, in order to study the inequality
$-Lu\ge f$, make the same assumptions on the operator $L$.
Namely, for a fixed $\eta\in\RN$ they say that {\bf (H$_\eta$)} is satisfied
if there exist a real number $Q=Q(\eta)>2$ and a nonnegative
continuous function $d_\eta:\RN\to\R_+$ such that the following
four properties hold
\begin{enumerate}
\item $d_\eta(\xi)=0$ if and only if $\xi=\eta$.
\item $d_\eta\in\Cdue(\RN\setminus\{\eta\})$.
\item The fundamental solution of $-L$ on $\RN$ at $\eta$
  is given by $\Gamma_\eta=d_\eta^{2-Q}$. That is, the functions
  $d_\eta^{2-Q}$ and $d_\eta^{1-Q}$ belong to $L^1_{loc}(\RN)$ and
  for any $\phi\in\Cdue_0(\RN)$ we have
  $$\int_{\RR^N}(-L\phi)(\xi)\Gamma_\eta(\xi)d\xi=\phi(\eta).$$
\item For any $i,j=1,\dots,l$  the functions $X_id_\eta$,
$X_j(d_\eta X_id_\eta)$ are bounded and
$|\grl d_\eta|^2\neq0$ almost everywhere on $\RN$.
\end{enumerate}
In this setting, it is immediate to check that the hypotheses of
Theorem \ref{th:hm} are fulfilled and a Hardy inequality related
to the operator $L$ holds.
\bt Let $\eta\in\RN$ and assume that  {\bf (H$_\eta$)} is satisfied.
    Then for any $u\in\Cuno_0(\RN)$, we have
$$\left(\frac{Q(\eta)-2}{2}\right)^2\int_{\RR^N}\frac{u^2(\xi)}{d_\eta^2(\xi)}
       {\abs{\grl d_\eta(\xi)}}^2 d\xi
     \le\int_{\RR^N} \abs{\grl u(\xi)}^{2} d\xi.    $$
\et

\medskip

As particular case of Theorem \ref{th:hm}, we obtain the following

\bt Let $g\in\Cdue(\Omega)$ be an $L_p$-harmonic function, that is $L_pg=0$
  and let $v\in\Cdue(\R)$ be a concave function such that
  $v\circ g$ is positive on $\Omega$.
  For any $u\in\Cuno_0(\Omega)$ we have
  $$\left(\frac{p-1}{p}\right)^p\int_\Omega\abs{u(\xi)}^p\frac{\abs{v'(g(\xi))}^p}{v^p(g(\xi))}
    \abs{\grl g(\xi)}^p\ d\xi\le\int_\Omega\abs{\grl u(\xi)}^p\ d\xi. $$
\et

The above result follows from Theorem \ref{th:hm} and the following worthwhile lemma
\bl\label{lem:harm} Let $p>1$, $\alpha\in\R$, $\alpha\neq 0$ and $g\in\Cdue(\Omega)$
  be a positive function such that $L_p(g^\alpha)=0$.
  Let $v\in\Cdue(\R)$ be such that $v'(g(\xi))\neq 0$ for $\xi\in\Omega$.
  Then setting $u(\xi)\decl v(g(\xi))$, we have
  $$L_pu=(p-1)\abs{\grl g}^p\abs{v'(g)}^{p-2}
   \left[ v''+\frac{1-\alpha}{g}v'\right].$$
\el
\bp{} We notice that for every $f\in\Cuno(\Omega)$ and $h\in\Cuno(\Omega,\R^l)$,
  we have $-\grl^*(fh)=\grl f\cdot h-f\grl^* h$.
  Let $\phi\in\Cuno(\R)$. Observing that $\grl (\phi\circ g)=\phi'(g) \grl g$,
  we deduce
  \begin{equation}
    \label{eq:tec1}
    L_p\phi=-\grl^*\left(\abs{\phi'(g)}^{p-2}\abs{\grl g}^{p-2}\phi'(g)\grl g\right)=
    \abs{\phi'(g)}^{p-2} \left[ (p-1) \phi'' \abs{\grl g}^p  +
      \phi' L_p g \right].
  \end{equation}
  Choosing $\phi(t)=t^\alpha$ in (\ref{eq:tec1}) and taking into account
  the $p$-harmonicity of $g^\alpha$
  we obtain $$gL_pg=(p-1)(1-\alpha)\abs{\grl g}^p,$$
  which substituted in (\ref{eq:tec1}), yields the claim.
\ep
\boss Lemma \ref{lem:harm} provides a generalization of the expression
of the usual $p$-Laplacian for radial function.
\eoss
\medskip

As simple application of previous results is the following
\begin{proposizione} Let $\grl$ be the usual gradient in $\R^2$, that is
  $\grl\decl\nabla = (\partial_x,\partial_y)^T$ and
  $\Omega\decl ]-\frac\pi2,\frac\pi2[\times\R$.
  For every $u\in\Cuno_0(\Omega)$ we have
  \be\label{dis:cos} \frac14 \int_\Omega \frac{u^2(x,y)}{\cos^2 x}\ dxdy
    \le\int_\Omega\abs{\nabla u(x,y)}^2\ dxdy.  \ee
  Moreover the constant $1/4$ is optimal and it is not achieved in $D^{1,2}(\Omega)$.
\end{proposizione}
\bp{} The inequality (\ref{dis:cos}) follows from Theorem \ref{th:hm} and
  the choice $d(x,y)\decl e^y\cos x$.

  Merging (\ref{dis:cos}) and the inequality
  $$\cos x\le \frac\pi2-\abs x=dist((x,y),\partial\Omega)\ \ \mathrm{for}\ \ (x,y)\in\Omega,$$
  we have
  \be\label{dis:cos2} \frac14 \int_\Omega \frac{u^2(x,y)}{dist^2((x,y),\partial\Omega)}\ dxdy
    \le\int_\Omega\abs{\nabla u(x,y)}^2\ dxdy.  \ee

   The fact that  $1/4$ is the best constant in (\ref{dis:cos2})
   (see \cite{mat-sob97}), implies the optimality of  $1/4$ in (\ref{dis:cos}).
   Hence applying Theorem \ref{th:achieve} we conclude the proof.
\ep

\boss
Finally, we notice that the result stated in our main Theorem \ref{th:hm}
can be reformulated also for non compact Riemannian manifold. This allow us
to re-obtain the Hardy inequalities present in \cite{car97} as well as their
many generalizations with weaker hypotheses.
\eoss

\section{Hardy Inequalities for some Subelliptic Operators}\label{sec:apply}
In this section we shall apply the previous results to
particular operators.

Let $1< p< \infty$. In the sequel for a given vector field $\grl$
and open set $\Omega\subset \RN$, we shall denote by
$D^{1,p}_L(\Omega)$ the closure of $\C_0^{\,\infty}(\Omega)$ in the
norm $(\int_\Omega\abs{\grl u}^p d\xi)^{1/p}$.
If $w\in L^1_{loc}(\Omega)$ and $w>0$ a.e. on $\Omega$,
$D^{1,p}_L(\Omega,w)$ denotes the closure of $\C_0^{\,\infty}(\Omega)$
in the norm $(\int_\Omega\abs{\grl u}^p w d\xi)^{1/p}$.

\subsection{Baouendi-Grushin operator}
Let $\RN$ be splitted in $\xi=(x,y)\in\R^n\times\R^k$. Let
$\gamma\ge0$ be a nonnegative real number and let $\mu$ be the matrix defined in
(\ref{matgrushin}). The corresponding vector field is
$\grgr=(\nabla_x,\abs x^\gamma\nabla_y)$ and the nonlinear
operator $L_p$ is $L_pu=\diverl(\abs{\grgr u}^{p-2}\grgr u)$. The
linear operator $L=L_2$ is the so-called Baouendi-Grushin   operator
$L=\lgr=\Delta_x+\abs x^{2\gamma}\Delta_y$. Notice that if $k=0$ or
$\gamma=0$, then $L$ and $L_p$ coincide  respectively
  with the usual Laplacian operator and $p$-Laplacian operator.

  Defining on $\R^N$ the dilation $\delta_\lambda$ as
  \be \delta_\lambda(x,y)\decl(\lambda x,\lambda^{1+\gamma}y);\label{dil:gr}\ee
  it is not difficult to check that
  $\grgr$ is homogeneous of degree one with respect
  to the dilation:
  $\grgr(\delta_\lambda)=\lambda\delta_\lambda(\grgr)$.

Let $\normg{\xi}=\normg{(x,y)}$ be the following distance from
the origin on $\R^N$:
$$\normg{\xi}=\normg{(x,y)}:=\left((\sum_{i=1}^dx_i^2)^{1+\gamma}+
    ({1+\gamma})^2\sum_{i=1}^ky_i^2\right)^{\frac{1}{2+2\gamma}}.$$
It is easy to see that $\normg{\cdot}$ is homogeneous of degree
one with respect to $\delta_\lambda$.

Let $Q\decl n+(1+\gamma)k$ be the so called  \emph{homogeneous dimension}.
If for $R>0$ we  denote by $B_R$ the set
\[B_{R}\decl\{\xi\in \R^N : \normg{\xi}<R\},\]
then we have
\[ \abs{B_R}=\abs{B_1} R^Q.\label{omog_mis}\]
Moreover, we have that  $\norg^s\in L^1_{loc}(\RN)$ if and only if $s>-Q$
and  $\norg^s\in L^1(\RN\setminus B_1)$ if and only if $s<-Q$.

The function $\normg\cdot$  is related to the fundamental solution at the
origin of Grushin operator $L$ (see \cite{dam-luc03}).
Namely, if $Q>2$ then the function $u_2\decl\norg^{2-Q}$ satisfies the relation
\[-L u_2=l_2\delta_0\quad\mathrm{on}\ \RN\]
in weak sense, where $\delta_0$ is the Dirac measure at $0$ and $l_2$
is a positive constant.

It is immediate to check that
\be \abs{\grgr \normg\xi}=\frac{|x|^\gamma}{\normg\xi^\gamma}\le 1.\label{grad_rad}\ee

Let  $p>1$ and let $\Gamma_p$ be the function defined as
\begin{equation}
  \Gamma_p(\xi)\decl\left\{
    \begin{array}{ll}
      \norg^\frac{p-Q}{p-1}& \mathrm{if}\ p\neq Q\\
      -\ln\norg &  \mathrm{if}\ p= Q
    \end{array} \right.\quad \mathrm{for}\ \xi \neq 0.
\end{equation}
A direct calculation gives that $\Gamma_p$ is $L_p$ harmonic on
$\RN\setminus\{0\}$,  that is
\begin{equation}
  \label{eq:gruhar}
  -L_p \Gamma_p = 0 \quad \mathrm{on}\ \RN\setminus\{0\}.
\end{equation}
Moreover, with analogous computations of
\cite{bal-tys02}, it is possible to show  that there exists a constant
$l_p\neq 0$ such that
\begin{equation}
  \label{eq:grufun}
  -L_p \Gamma_p=l_p\delta_0  \quad \mathrm{on}\ \RN
\end{equation}
in weak sense and  $l_p>0$ if and only if $Q\ge p>1$ (see also \cite{bie-gon04}).

\bt\label{teo:hargr} Let $\beta\in\R$ be fixed.
\begin{enumerate}
\item Let $1<p< +\infty$ and let $\Omega\subset\RN$ be an open set.
  If $\beta+Q<0$, then we also require that $0\not\in\Omega$.
  We have
  \begin{equation}\label{eq:hardyg}
   c_\beta^p\int_\Omega  {\abs{u}^p}{\normg\xi^{ \beta}}
                 \frac{|x|^{\gamma p}}{\norg^{\gamma p}}d\xi
                 \le\int_{\Omega} \abs{\grgr u}^{p}\norg^{p+\beta} d\xi
                 \quad\quad u\in D_L^{1,p}(\Omega,\norg^{p+\beta}),
  \end{equation}
  where $c_{\beta}\decl\frac{\abs{Q+\beta}}{p}$.
  In particular, we obtain
  \begin{equation} \left(\frac{\abs{Q-p}}{p}\right)^p\int_\Omega
    \frac{\abs{u}^p}{\normg\xi^{ p}}
                        \frac{|x|^{\gamma p}}{\normg\xi^{\gamma p}}d\xi
       \le \int_{\Omega} \abs{\grgr u}^{p} d\xi
       \quad\quad u\in D_L^{1,p}(\Omega) . \label{hardygpsi}
  \end{equation}
 Moreover, if $\Omega\cup \{0\}$ is a neighbourhood of the origin,
 then the constant $c_\beta^p$ is sharp.
\item  Let $p = Q>1$.
  Let $R>0$ and set $\Omega\decl\{\xi\in\RN|\ \norg<R\}$.
  If $\beta<-1$, then we have
 \begin{equation}\label{eq:hardygln}
   \tilde c_\beta^p\int_\Omega  \frac{\abs{u}^p}{\normg\xi^p}
       (\ln\frac R{\norg})^\beta \frac{|x|^{\gamma p}}{\norg^{\gamma p}}d\xi
     \le\int_{\Omega} \abs{\grgr u}^{p} (\ln\frac R{\norg})^{p+\beta} d\xi
     \quad\quad u\in D_L^{1,p}(\Omega,(\ln(R/\norg))^{p+\beta}),
  \end{equation}
  where $\tilde c_\beta\decl\frac{\abs{\beta+1}}{p}$.
  In particular, we obtain
 \begin{equation}\left(\frac{p-1}{p}\right)^p
   \int_\Omega  \frac{\abs{u}^p}{(\normg\xi\ln(R/\norg))^p}
       \frac{|x|^{\gamma p}}{\norg^{\gamma p}}d\xi
     \le\int_{\Omega} \abs{\grgr u}^{p} d\xi
     \quad\quad u\in D_L^{1,p}(\Omega).
  \end{equation}
 Moreover, the constant $\tilde c_\beta^p$ is sharp.
\end{enumerate}
\et

\boss If $\gamma=0$, then the operator $\Delta_\gamma$ is
  the standard Laplacian operator acting on functions defined on $\R^N$
  and \refe{hardygpsi} is the classical
  Hardy inequality (see (\ref{hardyeu}) where $w(\xi)=\abs{\xi}$).
\eoss
\boss The above inequality (\ref{eq:hardyg}) is already obtained in the case $1<p<Q$
  by the author in \cite{dam04-2}.
\eoss
\boss The inequality (\ref{eq:hardygln}) holds also for any $\beta>-1$
provided the set $\Omega$ is replaced by
$\Omega\decl\{\xi\in\RN|\ 0< \norg<R\}$.
\eoss

\bp{} We shall prove the inequalities for $u\in\Cuno_0(\Omega)$.
  The general case will follows by density argument.

  The inequalities (\ref{eq:hardyg}) for $p\neq Q$, and (\ref{eq:hardygln})
  follow from Theorem \ref{th:hm} choosing $d^\alpha=\Gamma_p$.

  Now, we prove the missing inequality (\ref{eq:hardyg}) when $p=Q$.
  We consider the case $Q+\beta>0$, the converse case is similar.
  Let $\sigma>0$ be such that $Q+\beta-\sigma>0$.
  We choose  $d(\xi)=\norg$ and $\alpha=\frac{\sigma}{Q-1}$.
  In this case it easy to check that $d^\alpha$ is sub-$L_Q$-harmonic
  on $\RN$, that is
  $$-L_Q(d^\alpha)=-\diverl \left((\frac{\sigma}{Q-1})^{Q-1}\norg^\sigma
    \frac{\abs{\grl \norg}^{Q-2} }{\norg^{Q-1}} \grl\norg \right)\le 0$$
  in weak sense.
  Indeed, using (\ref{eq:gruhar}) with $p=Q$, we deduce
  \begin{eqnarray*}
-L_Q(d^\alpha)&=&-\left(\frac{\sigma}{Q-1}\right)^{Q-1}\norg^\sigma
    \diverl \left(
        \Gamma_Q \right)
   - \sigma\left(\frac{\sigma}{Q-1}\right)^{Q-1} \frac{\abs{\grl \norg}^{Q} }{\norg^{Q-\sigma}} \\
   &=&\left(\frac{\sigma}{Q-1}\right)^{Q-1} \norg^\sigma\frac{\delta_0}{l_Q}-
   \frac{\sigma^Q}{{(Q-1)}^{Q-1}}
\frac{\abs{\grl \norg}^{Q} }{\norg^{Q-\sigma}}.
  \end{eqnarray*}
  Therefore, we get $-L_Q(d^\alpha)=-\frac{\sigma^Q}{{(Q-1)}^{Q-1}}
\frac{\abs{\grl \norg}^{Q} }{\norg^{Q-\sigma}}\le 0$ in weak sense.
  The constant $c$ in (\ref{eq:shar}) is
  $c=\frac{\sigma}{Q-1}\ \frac{\sigma-Q-\beta}{Q}<0$.

  These choices yield the inequality
$$(\frac{Q+\beta-\sigma}{Q})^p \int_\Omega  {\abs{u}^p}{\normg\xi^{ \beta}}
    \frac{|x|^{\gamma p}}{\norg^{\gamma p}}d\xi
    \le\int_{\Omega} \abs{\grgr u}^{p}\norg^{p+\beta} d\xi. $$
  Letting $\sigma\to 0$ we get the thesis.

  It remains to show that the constants $c_\beta^p$ and $\tilde c_\beta^p$
  appearing in (\ref{eq:hardyg}) and in (\ref{eq:hardygln}) are sharp.
  First we consider the case $\Omega=\RN\setminus\{0\}$.
  To this end it suffices to show that we are in the position to apply
  the scheme outlined in Remark \ref{rm:sharp}.

  Indeed, let $d(\xi)\decl\norg$ and for $\epsilon>0$ consider the constant
  $c(\epsilon)=\frac{Q+\beta}{p}+\frac{\epsilon}{p}$ defined in (\ref{eq:cepsilon}). We have
  $$c(\epsilon)p+\beta=\abs{Q+\beta}+Q+\beta+Q-Q+\epsilon  \ge -Q+\epsilon $$
and
  $$-c(\epsilon)p+\beta=-\abs{Q+\beta}+\beta+Q-Q-\epsilon \le -Q- \epsilon.$$
  These inequalities imply the integrability conditions
  (\ref{int:sharp}).
  Hence we obtain the sharpness of the constants.

  In order to conclude in the general case we proceed as follows:
  let $c_b(\Omega)$ be the best constant in \refe{eq:hardyg}.
By invariance of \refe{eq:hardyg} under the dilation
$\delta_\lambda$ defined in \refe{dil:gr} we have,
$$ c_b(B_R)=c_b(B_1)\ \ \mathrm{and}\ \ c_b(B_R\setminus\{0\})=c_b(B_1\setminus\{0\})\quad{\rm for\ any\ }R>0.$$
We note that if $B_R\setminus\{0\}\subset\Omega\subset \R^N\setminus\{0\}$ then,
\be c_{\beta}^p=c_b(\R^N\setminus\{0\})\le c_b(\Omega)\le c_b(B_R\setminus\{0\})=c_b(B_1\setminus\{0\}).\label{relbc0}\ee
Let $\phi\in\C^{\,\infty}_0(\R^N\setminus\{0\})$.
Since the support of $\phi$ is compact, then \refe{eq:hardyg} holds for $\phi$
with $\Omega=B_R\setminus\{0\}$, $R$ large enough and $c=c_b(B_R\setminus\{0\})=c_b(B_1\setminus\{0\})$.
Therefore  $c_b(B_1\setminus\{0\})\le c_b(\R^N\setminus\{0\})$, and from \refe{relbc0} we have $c_b(B_1\setminus\{0\})=c_\beta^p$.
Finally, since $B_R\subset\Omega\subset\RN$ imply
$$c_{\beta}^p\le c_b(\R^N)\le c_b(\Omega)\le c_b(B_R)=c_b(B_1)\le c_b(B_1\setminus\{0\})=c_\beta^p,$$
we conclude the proof.

  The optimality of the constant $\tilde c_\beta^p$ in
  (\ref{eq:hardygln}) can be easily proved using the procedure of
        Remark \ref{rm:sharp}.
\ep

Other inequalities of Hardy type related to Baouendi-Grushin
operator are given in the following
\bt\label{teo:hargr2} Let
$1\le m\le n$ and let $\beta\in\R$ be fixed.
  We set $z\decl(x_1,\dots,x_m)$.
\begin{enumerate}
\item  Let $1<p<+\infty$ and let $\Omega\subset\RN$ be an open set.
  If $m+\beta<0$, we also require that $\Omega\subset(\R^m\setminus\{0\})\times\R^{N-m}$.
  Then for every $u\in\Cuno_0(\Omega)$, we have
  \begin{equation}\label{eq:hardyg2}
   b_\beta^p\int_\Omega  {\abs{u}^p}{\abs z^{ \beta}} d\xi
                 \le\int_{\Omega} \abs{\grgr u}^{p}\abs z^{p+\beta} d\xi,
  \end{equation}
  where $b_{\beta}\decl\frac{\abs{m+\beta}}{p}$.

  In particular, for every $u\in\Cuno_0(\Omega)$, we obtain
  \begin{equation}
    \left(\frac{\abs{m-p}}{p}\right)^p\int_\Omega\frac{\abs{u}^p}{\normg\xi^p}d\xi
    \le\left(\frac{\abs{m-p}}{p}\right)^p\int_\Omega\frac{\abs{u}^p}{\abs z^p}d\xi
       \le \int_{\Omega} \abs{\grgr u}^{p} d\xi. \label{eq:hardygpsi2}
  \end{equation}
  Moreover, denoting with $B^s_r$ the Euclidean ball in $\R^s$ of radius $r$
  and center at the origin, if $B^m_r\setminus\{0\}\times B^{N-m}_r\subset \Omega$
  for some $r$, then the constant $b_\beta^p$ in (\ref{eq:hardyg2}) is sharp.
\item  Let $p = m>1$.
  Let $R>0$ and set $\Omega\decl\{(z,x_{m+1}\dots,x_n,y)\in\R^m\times\R^{n-m}\times \R^k,\ \abs  z<R\}$.
  If $\beta<-1$, then for every $u\in\Cuno_0(\Omega)$, we have
 \begin{equation}\label{eq:hardygln2}
   \tilde b_\beta^p\int_\Omega  \frac{\abs{u}^p}{\abs z^p}
       (\ln\frac R{\abs z})^\beta d\xi
     \le\int_{\Omega} \abs{\grgr u}^{p} (\ln\frac R{\abs z})^{p+\beta} d\xi,
  \end{equation}
  where $\tilde b_\beta\decl\frac{\abs{\beta+1}}{p}$.

  In particular, for every $u\in\Cuno_0(\Omega)$, we obtain
 \begin{equation}\left(\frac{p-1}{p}\right)^p
   \int_\Omega  \frac{\abs{u}^p}{(\abs z\ln(R/\abs z))^p}d\xi
     \le\int_{\Omega} \abs{\grgr u}^{p} d\xi.
  \end{equation}
  Moreover, the constant $\tilde b_\beta^p$ is sharp.
\end{enumerate}
\et

\bp{} The inequalities (\ref{eq:hardyg2})  and (\ref{eq:hardygln2})
  are a direct consequence of Theorem \ref{teo:harspec}.

  The fact that $\norg\ge\abs z$ yields the inequality (\ref{eq:hardygpsi2}).

  The sharpness of involved constants cannot be proved using
  the procedure of Remark \ref{rm:sharp}.
  Thus, we shall use a modification of the idea presented in
  \cite{sec-sme-wil03}.

  We prove the optimality of the constant $b_\beta^p$ in
  (\ref{eq:hardyg2}).
  The proof of the sharpness of the constant $\tilde b_\beta^p$ is similar.

  Let $c_b(\Omega)$ be the best constant in (\ref{eq:hardyg2}), that is
  \begin{equation}
    \label{eq:bestgru}
    c_b(\Omega)\decl \inf \left\{\frac{\int_\Omega\abs{\grgr \phi}^p\abs z^{\beta+p}}
       {\int_\Omega\abs \phi^p\abs z^\beta},\ \  \phi\in\Cuno_0(\Omega),\ \phi\neq 0\right\}.
  \end{equation}
  From (\ref{eq:hardyg2}) we have $c_b(\Omega)\ge b_\beta^p$.
  We shall prove the equality sign holds.

  First we consider the case
  $\Omega\decl\R^m\setminus\{0\}\times\R^{n-m}\times\R^{k}$.
  Observe that if we get the claim for $\Omega$, that is
  $c_b(\Omega)=b_\beta^p$, from $b_\beta^p\le c_b(\RN)\le c_b(\Omega)$,
  we get the claim also for $\Omega=\RN$.

  In what follows $t$ stands for the variables
  $t\decl(x_{m+1},\dots,x_n)\in\R^{n-m}$.
  Let $\phi\in\Cuno_0(\Omega)$ be such that $\phi=uvw$ with $u=u(z)$
  $v=v(t)$, $w=w(y)$, $v\in\Cuno_0(\R^m\setminus\{0\})$
  $v\in\Cuno_0(\R^{n-m})$ and  $w\in\Cuno_0(\R^k)$.
  It is clear that if $m=n$, then we choose $\phi=uw$ and the following
  proof results to be slightly simpler.

  By the convexity of the function $(q^2+r^2+s^2)^{p/2}$ for
  $q,r,s\ge 0$ we have
$$    (q^2+r^2+s^2)^{p/2}\le(1-\lambda-\mu)^{1-p}q^p+\lambda^{1-p}r^p+\mu^{1-p}s^p\quad \mathrm{for}\ \ \lambda,\mu>0,\ \lambda+\mu<1. $$
  Hence, for  $\lambda,\mu>0$ such that $\lambda+\mu<1$, we get
  \begin{eqnarray*}
    \lefteqn{\abs{\grgr \phi}^p=(v^2w^2\abs{\nabla_z u}^2+u^2w^2\abs{\nabla_t v}^2+u^2v^2\abs{\grgr w}^2)^{p/2}}\\
    &&\le (1-\lambda-\mu)^{1-p}\abs v^p\abs w^p\abs{\nabla_z u}^p
       +\lambda^{1-p}\abs u^p\abs w^p\abs{\nabla_t v}^p
       +\mu^{1-p}\abs u^p \abs v^p\abs{\grgr w}^p
\label{eq:tec2}
  \end{eqnarray*}
  Therefore, we obtain
  \begin{eqnarray*}
c_b(\Omega)\!\!&\le\!&\!\frac{\int_\Omega\abs{\grgr \phi}^p\abs z^{\beta+p}}
            {\int_\Omega\abs \phi^p\abs z^\beta}
     \le(1-\lambda-\mu)^{1-p}\frac{\int_\Omega \abs v^p\abs w^p\abs{\nabla_z u}^p\abs z^{\beta+p}}
            {\int_\Omega\abs v^p\abs w^p\abs u^p\abs z^\beta}\\
&&\quad+\lambda^{1-p}
            \frac{\int_\Omega \abs u^p\abs w^p\abs{\nabla_t v}^p\abs z^{\beta+p}}
            {\int_\Omega\abs u^p\abs v^p\abs w^p\abs z^\beta}+\mu^{1-p}
            \frac{\int_\Omega \abs u^p\abs v^p\abs{\grgr w}^p\abs z^{\beta+p}}
            {\int_\Omega\abs u^p\abs v^p\abs w^p\abs z^\beta} \\
&\!\!\le\!&\!(1-\lambda-\mu)^{1-p}\frac{\int_{\RR^m}\abs{\nabla_z u}^p\abs z^{\beta+p}dz}
              {\int_{\RR^m}\abs u^p\abs z^\beta dz}
              +\lambda^{1-p}
            \frac{\int_{\RR^{n-m}}\abs{\nabla_t v}^p dt}{\int_{\RR^{n-m}}\abs v^p dt}\
            \frac{ \int_{\RR^{m}}\abs u^p\abs z^{\gamma+\beta+p}dz}
               {\int_{\RR^{m}}\abs u^p\abs z^\beta dz}\\
&&\qquad\qquad +\mu^{1-p}
              \frac{\int_{\RR^k}\abs{\nabla_y w}^p dy}{\int_{\RR^k}\abs w^p dy}
              \frac{ \int_{\RR^n}\abs u^p \abs v^p\abs z^{\beta+p}\abs x^{\gamma p}dx}
               {\int_{\RR^n}\abs u^p\abs v^p\abs z^\beta dx}
  \end{eqnarray*} 
  Now, the infimum of the ratio
  ${\int_{\RR^k}\abs{\nabla_y w}^p dy}/\int_{\RR^k}\abs w^p dy$
  vanishes, as well as the ratio
  ${\int_{\RR^{n-m}}\abs{\nabla_t v}^p dt}/\int_{\RR^{n-m}}\abs v^p dt$.
  From the classical Hardy inequalities (see also Theorem \ref{teo:hargr} with $k=0$, $n=N=Q$), the infimum of
  ${\int_{\RR^m}\abs{\nabla_z u}^p\abs z^{\beta+p}dz}/
              {\int_{\RR^m}\abs u^p\abs z^\beta dz}$ is $b_\beta^p$.
  Thus, letting $\lambda,\mu\rightarrow 0$, we get the claim.

  In order to complete the proof, we prove the claim in the case
  $B^m_r\setminus\{0\}\times B^{N-m}_r\subset \Omega\subset\R^m\setminus\{0\}\times \R^{N-m}$ for some $r>0$.

  Let $B^*_r\decl B^m_r\setminus\{0\}\times B^{n-m}_r\times
    B^k_{r^{1+\gamma}}$.
  For $s>0$ sufficiently small we have that
  $B^*_s=B^m_s\setminus\{0\}\times B^{n-m}_s\times B^k_{s^{1+\gamma}}\subset\Omega\subset\R^m\setminus\{0\}\times \R^{n-m}$.
  Thus, we obtain $b_\beta^p=c_b(\R^m\setminus\{0\}\times \R^{n-m})
  \le c_b(\Omega)\le c_b(B^*_s)$.
  By invariance of \refe{eq:hardyg2} under the dilation
  $\delta_\lambda$ defined in  \refe{dil:gr} we have,
  $c_b(B^*_r)=c_b(B^*_1)\quad{\rm for\ any\ }r>0$.
  Arguing as in the proof of Theorem \ref{teo:hargr},
  we get the claim and conclude the proof.
\ep


\subsection{Heisenberg-Greiner operator}
\label{sec:hei-gre}
   Let $\xi=(x,y,t)\in\R^n\times\R^n\times\R$, $r\decl\abs{(x,y)}$,
   $\gamma\ge1$ and let $\mu$ be the matrix defined in
   (\ref{matgrein}).
  We remind that for $p=2$ and $\gamma=1$ $L_p$ is the sub-Laplacian $\lh$ on
  the Heisenberg group $\h^n$.
  If $p=2$ and $\gamma=2,3,\dots$, $L_p$ is a Greiner operator
  (see \cite{gre79}).

For $(x,y,t)\in\R^n\times\R^n\times\R$,
we define $$N(x,y,t)\decl((x^2+y^2)^{2\gamma}+t^2)^{1/4\gamma}=(r^{4\gamma}+t^2)^{1/4\gamma},$$
where we have set $r\decl(x^2+y^2)^{1/2}$.
Let $Q\decl 2n+2\gamma$, $p>1$ and let $\Gamma_p$ be the function defined as
\begin{equation}
  \Gamma_p\decl\left\{
    \begin{array}{ll}
      N^\frac{p-Q}{p-1}& \mathrm{if}\ p\neq Q\\
      -\ln N &  \mathrm{if}\ p= Q
    \end{array} \right.\quad \mathrm{for}\ \xi \neq 0.
\end{equation}
The function $\Gamma_p$ is $L_p$ harmonic on
$\RN\setminus\{0\}$,  that is
\begin{equation}
  -L_p \Gamma_p = 0 \quad \mathrm{on}\ \RN\setminus\{0\}.
\end{equation}
Moreover, arguing as in \cite{bal-tys02}, there exists a constant $l_p\neq 0$ such that
\begin{equation}
  \label{eq:greinfun}
  -L_p \Gamma_p=l_p\delta_0
\end{equation}
in weak sense and $l_p>0$ if and only if  $Q\ge p$ (see also
\cite{zha-niu03}). Moreover, $\abs{\grl
N}=\frac{r^{2\gamma-1}}{N^{2\gamma-1}}$.

\bt\label{teo:hargrein} Let  $\beta\in\R$ be fixed.
\begin{enumerate}
\item Let $1<p<+\infty$ and let $\Omega\subset\RN$ be an open set.
  If $\beta+Q<0$ we also require that $0\not\in\Omega$.
  Then, we have
  \begin{equation}\label{eq:hardygrein}
   c_\beta^p\int_\Omega  {\abs{u}^p}{N^{ \beta}}
                 \frac{r^{p(2\gamma-1)}}{N^{p(2\gamma-1)}}d\xi
                 \le\int_{\Omega} \abs{\grl u}^{p}N^{p+\beta} d\xi
                 \quad\quad u\in D_L^{1,p}(\Omega,N^{p+\beta}),
  \end{equation}
  where $c_{\beta}\decl\frac{\abs{Q+\beta}}{p}$.

  In particular, we obtain
  \begin{equation} \left(\frac{\abs{Q-p}}{p}\right)^p\int_\Omega
    \frac{\abs{u}^p}{N^{ p}}
                        \frac{r^{p(2\gamma-1)}}{N^{p(2\gamma-1)}}d\xi
       \le \int_{\Omega} \abs{\grl u}^{p} d\xi
       \quad\quad u\in D_L^{1,p}(\Omega). \label{hardygreinpsi}
  \end{equation}
 Moreover, if $\Omega\cup \{0\}$ is a neighbourhood of the origin,
 then the constant $c_\beta^p$ is sharp.
\item  Let $p = Q$.
  Let $R>0$ and set $\Omega\decl\{\xi\in\RN |\ N(\xi)<R\}$.
  If $\beta<-1$, then we have
 \begin{equation}\label{eq:hardygreinln}
   \tilde c_\beta^p\int_\Omega  \frac{\abs{u}^p}{N^p}
       (\ln\frac R{N})^\beta \frac{r^{p(2\gamma-1)}}{N^{p(2\gamma-1)}}d\xi
     \le\int_{\Omega} \abs{\grl u}^{p} (\ln\frac R{N})^{p+\beta} d\xi
     \quad u\in D_L^{1,p}(\Omega,(\ln(R/N))^{p+\beta}),
  \end{equation}
  where $\tilde c_\beta\decl\frac{\abs{\beta+1}}{p}$.

  In particular, we obtain
 \begin{equation}\left(\frac{p-1}{p}\right)^p
   \int_\Omega  \frac{\abs{u}^p}{(N\ln(R/N))^p}
       \frac{r^{p(2\gamma-1)}}{N^{p(2\gamma-1)}}d\xi
     \le\int_{\Omega} \abs{\grl u}^{p} d\xi
     \quad\quad u\in D_L^{1,p}(\Omega).
  \end{equation}
 Moreover, the constant  $\tilde c_\beta^p$ is sharp.
\end{enumerate}
\et

\boss If $\gamma=1$, then the operator $L_p$ is the counterpart of the
  $p$-Laplacian for
  the sub-Laplacian operator
  acting on functions defined on the Heisenberg group $\h^n$.
  In this case the Hardy inequality (\ref{hardygreinpsi}) is already obtained
  for $1<p<Q$ by
  Garofalo and Lanconelli in \cite{gar-lan90},
  Niu, Zhang and Wang in \cite{niu-zha-wan01}.
  The author in \cite{dam04-1}
  proves the inequality (\ref{eq:hardygrein}) and the sharpness of
  the involved constant.

  In the general case $\gamma\ge 1$, the inequality (\ref{hardygreinpsi}) is
  already obtained in the case $1<p<Q$
  for function $u\in\Cuno_0(\RN\setminus \{0\})$ in \cite{zha-niu03}.
\eoss

The proof of the above theorem follows arguing as in the proof of
Theorem \ref{teo:hargr}.
Arguing as in Theorem \ref{teo:hargr2} we obtain the following
\bt\label{teo:hargrein2} Let $\beta\in\R$ be fixed.
\begin{enumerate}
\item Let $1<p<+\infty$ and let $\Omega\subset\RN$ be an open set.
  If $2n+\beta<0$, we also require that $\Omega\subset(\R^{2n}\setminus\{0\})\times\R$.
  Then for every $u\in\Cuno_0(\Omega)$, we have
  \begin{equation}\label{eq:hardygrein2}
   b_\beta^p\int_\Omega  {\abs{u}^p}{ r^{ \beta}} d\xi
                 \le\int_{\Omega} \abs{\grl u}^{p}r^{p+\beta} d\xi,
  \end{equation}
  where $b_{\beta}\decl\frac{\abs{2n+\beta}}{p}$.
  Moreover, denoting with $B^s_r$ the Euclidean ball in $\R^s$ of radius $r$
  with center at the origin, if $B^{2n}_r\setminus\{0\}\times B^{1}_r\subset \Omega$
  for some $r$, then the constants $b_\beta^p$ is sharp.

  In particular, for every $u\in\Cuno_0(\Omega)$, we obtain
  \begin{equation}
    \left(\frac{\abs{2n-p}}{p}\right)^p\int_\Omega\frac{\abs{u}^p}{N^p}d\xi
    \le\left(\frac{\abs{2n-p}}{p}\right)^p\int_\Omega\frac{\abs{u}^p}{r^p}d\xi
       \le \int_{\Omega} \abs{\grl u}^{p} d\xi. \label{eq:hardygreinpsi2}
  \end{equation}
\item  Let $p =2n$.
  Let $R>0$ and set $\Omega\decl\{(x,y,t)\in\R^{2n}\times \R,\ |(x,y)|<R\}$.
  If $\beta<-1$, then for every $u\in\Cuno_0(\Omega)$, we have
 \begin{equation}\label{eq:hardygreinln2}
   \tilde b_\beta^p\int_\Omega  \frac{\abs{u}^p}{r^p}
       (\ln\frac R{r})^\beta d\xi
     \le\int_{\Omega} \abs{\grl u}^{p} (\ln\frac R{r})^{p+\beta} d\xi,
  \end{equation}
  where $\tilde b_\beta\decl\frac{\abs{\beta+1}}{p}$.

  In particular, for every $u\in\Cuno_0(\Omega)$, we obtain
 \begin{equation}\left(\frac{p-1}{p}\right)^p
   \int_\Omega  \frac{\abs{u}^p}{(\abs x\ln(R/r))^p}d\xi
     \le\int_{\Omega} \abs{\grl u}^{p} d\xi.
  \end{equation}
 Moreover, the constant $\tilde b_\beta^p$  is sharp.
\end{enumerate}
\et


\subsection{Hardy Inequalities on Carnot Groups}\label{sec:carnot}
In this section we shall present some Hardy inequalities
in the framework of Carnot Groups.

We begin by  quoting some preliminary facts on these
structures and refer the interested reader to \cite{bon-ugu05, fol75, fol-ste82, hei95})
for  more precise information  on this subject.

A Carnot group is a connected, simply connected, nilpotent Lie
group $\G$ of dimension $N$ with graded Lie algebra ${\cal
G}=V_1\oplus \dots \oplus V_r$ such that $[V_1,V_i]=V_{i+1}$ for
$i=1\dots r-1$ and $[V_1,V_r]=0$. A such integer $r$ is called the
\emph{step} of the group.
 We set $l=n_1=\dim V_1$, $n_2=\dim V_2,\dots,n_r=\dim V_r$.
A  Carnot group $\G$ of dimension $N$ can be identified, up to an
isomorphism, with the structure of a \emph{homogeneous Carnot
Group} $(\RN,\circ,\delta_\lambda)$ defined as follows; we
identify $\G$ with $\RN$ endowed with a Lie group law $\circ$. We
consider $\RN$ splitted in $r$ subspaces
$\RN=\R^{n_1}\times\R^{n_2}\times\cdots\times\R^{n_r}$ with
$n_1+n_2+\cdots+n_r=N$ and $\xi=(\xi^{(1)},\dots,\xi^{(r)})$ with
$\xi^{(i)}\in\R^{n_i}$. We shall assume that there exists a family
of Lie group automorphisms, called \emph{dilation},
$\delta_\lambda$ with $\lambda>0$ of the form
$\delta_\lambda(\xi)=(\lambda\xi^{(1)},\lambda^2
\xi^{(2)},\dots,\lambda^r \xi^{(r)})$. The Lie algebra of
left-invariant vector fields on $(\RN,\circ)$ is $\cal G$. For
$i=1,\dots,n_1=l $ let $X_i$ be the unique vector field in $\cal
G$ that coincides with $\partial/\partial\xi^{(1)}_i$ at the
origin. We require that the Lie algebra generated by
$X_1,\dots,X_{l}$ is the whole $\cal G$.

If the above hypotheses are satisfied, we shall call
$\G=(\RN,\circ,\delta_\lambda)$ a \emph{homogeneous Carnot Group}.
We denote with $\grl$ the vector field $\grl\decl(X_1,\dots,X_l)^T$.
The \emph{canonical sub-Laplacian} on $\G$ is the
second order differential operator defined by
$L_2\decl\Delta_G=\sum_{i=1}^{l} X_i^2$ and we define for $p>1$ the
$p$-sub-Laplacian operator
$L_p(u)\decl\sum_{i=1}^{l} X_i(\abs{\grl u}^{p-2}X_iu)$.

Some important  properties of Homogeneous Carnot groups are the
following: the Lebesgue measure on $\RN$ coincides  with the
bi-invariant Haar measure on {\doppio G}. We denote by
$Q\decl \sum_{i=1}^r i\,n_i=  \sum_{i=1}^r i\,\mathrm{dim}V_i$ the \emph{homogeneous dimension} of $\G$.
For every measurable set $E\subset\RN$, we have
$|\delta_\lambda(E)|=\lambda^Q|E|$. Since $X_1,\dots,X_{l}$
generate the whole $\cal G$, the sub-Laplacian $L$ satisfies the
H\"ormander's hypoellipticity  condition. Moreover, the vector
fields $X_1,\dots,X_{l}$ are homogeneous of degree 1 with respect
to $\delta_\lambda$.

A nonnegative continuous function $N:\RN\to\R_+$ is called a 
\emph{homogeneous norm} on {\doppio G}, if 
$N(\xi^{-1})=N(\xi)$, $N(\xi)=0$ if and only if $\xi=0$ and it is
homogeneous of degree 1 with respect to $\delta_\lambda$ (i.e.
$N(\delta_\lambda(\xi))=\lambda N(\xi)$).
A homogeneous norm $N$ defines on $\G$ a \emph{pseudo-distance} as
$d(\xi,\eta)\decl N(\xi^{-1}\eta)$. 
For such a function $d$, there holds only a pseudo-triangular inequality:
\be d(\xi,\eta)\le C d(\xi,\zeta)+C d(\zeta,\eta)\qquad(\xi,\zeta,\eta\in\G)\label{dis:pseudotr}\ee
with $C\ge 1$. Hence, $d$, in general, is not a distance.

If $N$ and $\tilde N$ are two homogeneous norms, then they are equivalent,
that is, there exists a constant
$C>0$ such that $C^{-1}N(\xi)\le \tilde N(\xi)\le CN(\xi)$.

Let $N$ be a homogeneous norm, then there exists a constant
$C>0$ such that $C^{-1}\abs\xi\le N(\xi)\le C\abs\xi^{1/r}$,
for $N(\xi)\le1$ and $\abs\cdot$ stands for the Euclidean norm. An
example of homogeneous norm is the following
\be
    N_S(\xi)\decl\left(\sum_{i=1}^r\abs{\xi_i}^{2r!/i}\right)^{1/2r!}.\label{eq:norms}\ee

Notice that if $N$ is a homogeneous norm differentiable a.e., 
then $\abs{\grl N}$ is homogeneous of degree 0 with respect to 
$\delta_\lambda$, hence $\abs{\grl N}$ is bounded.

\bigskip

Special examples of Carnot groups are the
  Euclidean spaces $\R^Q$.
  Moreover, if $Q\le 3$ then any Carnot group is the ordinary Euclidean
  space $\R^Q$.

 The most simple nontrivial example of a Carnot group
  is the Heisenberg group $\hei^1=\R^3$.
  For an integer $n\ge1$, the Heisenberg group $\hei^n$ is defined as follows:
  let $\xi=(\xi^{(1)},\xi^{(2)})$ with
  $\xi^{(1)}\decl(x_1,\dots,x_n,y_1,\dots,y_n)$ and $\xi^{(2)}\decl t$.
  We endow $\R^{2n+1}$ with  the group law
\[\hat\xi\circ\tilde\xi\decl(\hat x+\tilde x,\hat y+\tilde y,\hat t+ \tilde t+2\sum_{i=1}^n(\tilde x_i\hat y_i-\hat x_i \tilde y_i)).\]
For $i=1,\dots,n$, consider the vector fields
\[X_i\decl\frac{\partial}{\partial x_i}+2y_i\frac{\partial}{\partial t},\
        Y_i\decl\frac{\partial}{\partial y_i}-2x_i\frac{\partial}{\partial t},\]
and the associated  Heisenberg gradient as follows
\[ \grh\decl (X_1,\dots,X_n,Y_1,\dots,Y_n)^T.\]
The sub-Laplacian $\lh$ is then the operator defined by
\[ \lh\decl\sum_{i=1}^nX_i^2+Y_i^2.\]
The family of dilation is given by
\[\delta_\lambda(\xi)\decl (\lambda x,\lambda y,\lambda^2 t).\]
In ${\hei^n}$ we can define  the canonical homogeneous norm by
\[\abs{\xi}_H\decl \left(\left(\sum_{i=1}^n x_i^2+y_i^2\right)^2+t^2\right)^{1/4}.\]
The homogeneous dimension is given by $Q=2n+2$, and the
fundamental solution of the sub-Laplacian $-\lh$ at point $\eta$ takes the form
$\Gamma_\eta(\xi)=\abs{\eta^{-1}\circ \xi}_H^{-2n}$.

\bigskip

Other particular cases of Carnot groups are the Heisenberg-type groups.
They were introduced by Kaplan \cite{kap80} and have subsequently studied
by several authors. We list some properties for Heisenberg-type groups
and refer the reader to \cite{bon-ugu04, gar-vas01} and the reference therein.

Let $\G$ be a Carnot group of step 2 with Lie algebra
${\cal G}=V_1\oplus V_2$ and let $V_1$ be endowed with a scalar product
 $\langle\cdot,\cdot\rangle$.
Let $J:V_2\to End(V_1)$ be defined as
$$\langle J(\eta)\hat\xi,\tilde\xi\rangle =\langle [\hat\xi,\tilde\xi],\eta\rangle, \qquad
  \eta\in V_2,\quad \hat\xi,\ \tilde\xi\in V_1.$$

We say that $\G$ is of \emph{H(eisenberg)-type} if for all $\eta\in V_2$
we have $J(\eta)^2=-\abs \eta^2 Id$.

Let $\G$ be an $H$-type group. Denoting by $\exp$ the exponential map
$\exp:{\cal G}\to\G$ (that is a global diffeomorphism), we define
the analytic mappings $x:\G\to V_1$ and $t:\G\to V_2$ by the identity
$\xi=\exp(x(\xi)+t(\xi))$. For the sake of simplicity we shall identify
$\xi$ with $\xi=(x,t)$.

Let $N$ be defined as
\be \label{eq:normH}N(\xi)\decl(\abs{x}^4+16\abs{t}^2)^{1/4}.\ee
Then $N$ is a homogeneous norm on $\G$.
In this setting, the homogeneous dimension is given by $Q=n_1+2n_2$
(we remind that $n_1=\dim V_1$ and $n_2=\dim V_2$).

Let  $p>1$ and let $\Gamma_p$ be the function defined as
\begin{equation}
  \Gamma_p(\xi)\decl\left\{
    \begin{array}{ll}
      N^\frac{p-Q}{p-1}& \mathrm{if}\ p\neq Q\\
      -\ln N &  \mathrm{if}\ p= Q
    \end{array} \right.\quad \mathrm{for}\ \xi \neq 0.
\end{equation}
The function $\Gamma_p$ is $L_p$ harmonic on
$\RN\setminus\{0\}$,  that is
\begin{equation}
  -L_p \Gamma_p = 0 \quad \mathrm{on}\ \RN\setminus\{0\}.
\end{equation}
Moreover, there exists a constant $l_p$ such that
\begin{equation}
  \label{eq:htypefun}
  -L_p \Gamma_p=l_p\delta_0
\end{equation}
and if $Q\ge p$, then $l_p>0$ (see
\cite{cap-dan-gar96,hei-hol97, kap80}). Moreover, $\abs{\grl
N(\xi)}=\frac{\abs{x}}{N(\xi)}$.

Suppose that a function $u$ has the form $u=u(\abs x, t)$, then
we have
$$\abs{\grl u(\xi)}^2=\abs{\nabla_x u}^2 +
\frac{\abs{x}^2}{4}\abs{\nabla_t u}^2\quad\mathrm{ and}\quad
Lu(\xi)=\Delta_xu+\frac{\abs{x}^2}{4}\Delta_t u.$$

\bigskip

Now we come back to the general Carnot group.
It is well-known that there exists a homogeneous
norm $N_2$ smooth on $\G\setminus\{0\}$ such that
$(N_2(\xi))^{2-Q}$
is a fundamental solution of $-L_2$ at $0$
(see \cite{fol75, gal82}).
On the other hand there exists a homogeneous
norm $N_Q$ on $\G$ such that
$-\ln N_Q$
is a fundamental solution of $-L_Q$ at $0$
(see \cite{bal-hol-tys02,hei-hol97}).
In general these two norms do not agree (see \cite{bal-tys02}).
Moreover, according to author's knowledge, the best result on the regularity of
$N_Q$ is that it is H\"older continuous (\cite{cap-dan-gar93, cap-dan-gar96}, 
see also \cite{cap99}).

In spite of lack of information on regularity of $N_Q$,
we can still use the results of previous section to obtain
Hardy inequalities related to $\grl$ involving the homogeneous
norm $N_2$ and $N_Q$ for $p=2$ and $p=Q$.

In the case $1<p<Q$ one can argue as follows. Assume that $G_p$ is a 
fundamental solution of $-L_p$ at $0$ on $\G$ (that is $-L_pG_p=\delta_0$)
with a singularity at $0$. We set $N_p\decl G_p^{\frac{p-1}{p-Q}}$.
Now applying the results of previous section we get a Hardy inequality
involving the function $N_p$.
Using the results presented in \cite{cap-dan-gar96}, 
it is easy to prove that if $N$ is a homogeneous norm on $\G$,
then there exists a constant $C>0$ such that
$$ CN(\xi)\le N_p(\xi)\le C^{-1}N(\xi)\quad \mathrm{for\ every\ } \xi\in\G.$$
Hence, we obtain a Hardy inequality involving a homogeneous norm $N$, more precisely
\begin{proposizione} Under the above hypotheses, there exists a constant $c>0$
such that for every $u\in\Cuno_0(\G)$, we have
\be c\int_{\GG} \frac{\abs{u}^p}{N^{ p}} \abs{\grl N_p}^p d\xi
 \le \int_{\GG} \abs{\grl u}^{p} d\xi.\label{dis:harcarnp}\ee
\end{proposizione}

If, in the previous inequality (\ref{dis:harcarnp}) we fix, for instance, $N=N_2$,
we cannot say anything on the constant $c$ and, in particular,
we are not able to estimate  $c$:
This is due to the  lack of information about
the relation between $N_p$ and $N_2$.

Therefore, in what follows, for $p>1$ we denote with $\Gamma_p$ the function defined as
\begin{equation}
  \Gamma_p(\xi)\decl\left\{
    \begin{array}{ll}
      N_2^\frac{p-Q}{p-1}& \mathrm{if}\ p\neq Q\\
      -\ln N_2 &  \mathrm{if}\ p= Q
    \end{array} \right.\quad \mathrm{for}\ \xi \neq 0.\label{polfun}
\end{equation}

The question if $\Gamma_p$ is $L_p$ harmonic on $\G\setminus\{0\}$ arises.

In \cite{bal-tys02} the authors give the following definition
\bd The group $\G$ is \emph{polarizable} if
$N_2$ is $\infty$-harmonic on $\G\setminus\{0\}$, that is, $N_2$
is a solution of
$$\Delta_\infty f\decl \frac 12 \langle \grl\abs{\grl f}^2,\grl f\rangle=0
  \quad\mathrm{on}\quad \G\setminus\{0\}. $$
\ed
We recall that for $f\in\Cdue$ we can write $\Delta_\infty f$ also
as
$$\Delta_\infty f =  \langle (\grl^{2*} f)\grl f,\grl f\rangle$$
where $\grl^{2*}f$ denotes the symmetrized horizontal Hessian matrix of $f$,
$\grl^{2*}f\decl 1/2 [(\grl^{2}f) +(\grl^2f)^T]$.

In \cite{bal-tys02} the authors prove that if $\G$ is polarizable
then $\Gamma_p$ defined in (\ref{polfun})
is $p$-harmonic on $\G\setminus\{0\}$. Moreover, there exists
$l_p\neq 0$ such that $-L_p(\Gamma_p)=l_p\delta_0$ on $\G$ and
$l_p>0$ if and only if $1<p\le Q$. 

Actually, the condition that $\Gamma_p$ is $L_p$-harmonic on 
 $\G\setminus\{0\}$ is also a sufficient condition for the polarizability
as specified by the following
\begin{proposizione}\label{pr:pol} The group $\G$ is polarizable if and only if the function $\Gamma_p$
  defined in (\ref{polfun}) is $L_p$-harmonic on $\G\setminus\{0\}$ for some
  $p>1$, $p\neq 2$ (and hence for all $p>1$).
\end{proposizione}
\bp{} The necessary condition is already proved in \cite{bal-tys02}.
Thus, we shall prove the sufficient condition.

Let $u$ be a smooth function.  By computation we have
  \be\label{eq:lpexplic}
L_p u=\grl(\abs{\grl u}^{p-2})\cdot\grl u+
   \abs{\grl u}^{p-2}L_2 u=(p-2)\abs{\grl u}^{p-4}\Delta_\infty u
     + \abs{\grl u}^{p-2}L_2u. \ee
  Taking into account that $N_2^{2-Q}$ is $L_2$-harmonic on $\G\setminus\{0\}$,
  from Lemma \ref{lem:harm}, we have $L_2N_2=(Q-1)\frac{\abs{\grl N_2}^2}{N_2}$.
  Hence, applying (\ref{eq:lpexplic}) to $N_2$ we have
  \be L_p N_2=(p-2)\abs{\grl N_2}^{p-4}\Delta_\infty N_2+(Q-1)\frac{\abs{\grl N_2}^p}{N_2}.\label{id:lpn2}\ee
  The thesis will follow if we prove that the identity
  \be\label{id:goal1}L_pN_2=(Q-1)\frac{\abs{\grl N_2}^p}{N_2}
\ee
holds for every $\xi\neq 0$.

  Let $p>1$, $p\neq 2$ be such that  $\Gamma_p$ is $L_p$-harmonic  $\G\setminus\{0\}$.
  First we assume that $p\neq Q$.
  We apply Lemma \ref{lem:harm} with $g=N_2$, $\alpha=\frac{p-Q}{p-1}$
  to $u=N_2$ obtaining the identity (\ref{id:goal1}).

  Now we consider the case $p=Q$. Since $\Gamma_Q=-\ln N_2$ is $Q$-harmonic,
  the function $-\ln \frac{N_2}{R}$ is still  $Q$-harmonic and positive on
  $\Omega_R\decl\{\xi\in\G |\,0<N_2(\xi)<R\}$. 
  Thus applying Lemma \ref{lem:harm} with $g=-\ln \frac{N_2}{R}$, $\alpha=1$
  to $u=N_2$ we have that the identity (\ref{id:goal1}) is fulfilled on $\Omega_R$. 
  Since $R$ is arbitrary
  we conclude that the identity (\ref{id:goal1}) holds on $\G\setminus\{0\}$.
\ep

\medskip
Examples of polarizable Carnot groups are the usual Euclidean space,
as well as H-type group and hence the Heisenberg group.
This is proved in \cite{bal-tys02}.
\boss Proposition \ref{pr:pol} provides a straightforward
  proof of the polarizability of H-type groups.
\eoss

\bt\label{teo:harcarn} Let $p>1$ and 
let $\Gamma_p$ be $L_p$-harmonic on $\G\setminus\{0\}$.
Let $\beta\in\R$ be fixed and  let $N=N_2$.
\begin{enumerate}
\item Let $1<p< +\infty$ and let $\Omega\subset\G$ be an open set.
  If $\beta+Q<0$ we also require that $0\not\in\Omega$.
  Then  we have
  \begin{equation}\label{eq:hardycarn}
   c_\beta^p\int_\Omega  {\abs{u}^p}{N^{ \beta}}\abs{\grl N}^p d\xi
                 \le\int_{\Omega} \abs{\grl u}^{p}N^{p+\beta} d\xi
                 \quad\quad u\in D_L^{1,p}(\Omega,N^{p+\beta}),
  \end{equation}
  where $c_{\beta}\decl\frac{\abs{Q+\beta}}{p}$.
  In particular, we obtain
  \begin{equation} \left(\frac{\abs{Q-p}}{p}\right)^p\int_\Omega
    \frac{\abs{u}^p}{N^{ p}} \abs{\grl N}^p d\xi
       \le \int_{\Omega} \abs{\grl u}^{p} d\xi
       \quad\quad u\in D_L^{1,p}(\Omega). \label{hardycarnpsi}
  \end{equation}
  Moreover, if $\Omega\cup \{0\}$ is a neighbourhood of the origin,
  then the constant $c_\beta^p$ is sharp.
\item  Let $p = Q>1$.
  Let $R>0$ and set $\Omega\decl\{\xi\in\G,\ N(\xi)<R\}$.
  If $\beta<-1$, then  we have
 \begin{equation}\label{eq:hardycarnln}
   \tilde c_\beta^p\int_\Omega  \frac{\abs{u}^p}{N^p}
       (\ln\frac R N)^\beta\abs{\grl N}^p d\xi
     \le\int_{\Omega} \abs{\grl u}^{p} (\ln\frac RN)^{p+\beta} d\xi
     \quad\quad u\in D_L^{1,p}(\Omega,(\ln(R/N))^{p+\beta}),
  \end{equation}
  where $\tilde c_\beta\decl\frac{\abs{\beta+1}}{p}$.

  In particular, we obtain
 \begin{equation}\left(\frac{p-1}{p}\right)^p
   \int_\Omega  \frac{\abs{u}^p}{(N\ln(R/N))^p}\abs{\grl N}^p d\xi
     \le\int_{\Omega} \abs{\grl u}^{p} d\xi
     \quad\quad u\in D_L^{1,p}(\Omega).
  \end{equation}
  Moreover, the constant $\tilde c_\beta^p$ is sharp.
\end{enumerate}
\et

\boss The above theorem still holds for $p=2$ with $N=N_2$ and for $p=Q$
  with $N=N_Q$  in any Carnot Group and without the hypothesis
  of polarizability.
\eoss

If $\G=\RN$ and $\grl=\nabla$ is the usual gradient, then $\abs{\nabla N}=1$
and the above inequalities are a generalization of the known
Hardy inequalities.

\bigskip

Let $d_i=n_1+\dots+n_i$ for $i=1,\dots,r$ so that  $d_1=n_1=l$ and $d_r= N$.
It results
$$X_i=\derivpar{\xi_i}+\sum_{k=1}^{r-1} \sum_{s=d_k+1}^{d_{k+1}}
  P_{i,s,k}(\xi_1,\dots,\xi_{d_1},\xi_{d_1+1},\dots,\xi_{d_2},\dots,
    \xi_{d_{k-1}+1},\dots,\xi_{d_k})\derivpar{\xi_s},$$
where $P_{i,s,k}$ is a polynomial homogeneous of degree $k$ with
respect to dilation $\delta_\lambda$.
Denoting with $\mu$ the
matrix such that
$X_i=\sum_{j=1}^{N}\mu_{ij}(\xi)\frac{\partial}{\partial\xi_j}$,
it results that $\mu$ has the form $\mu=(I_l,\mu_1)$, hence
in particular $\mu$ has the form
(\ref{matspec}). Therefore we have the following
\bt\label{teo:harcarn2}Let $1\le m\le l$ and let $\beta\in\R$
be fixed. We set $z\decl(x_1,\dots,x_m)$.
\begin{enumerate}
\item  Let $1<p<+\infty$ and let $\Omega\subset\G$ be an open set.
  If $m+\beta<0$, we also require that $\Omega\subset(\R^m\setminus\{0\})\times\R^{N-m}$.
  Then for every $u\in\Cuno_0(\Omega)$, we have
  \begin{equation}\label{eq:hardycarn2}
   b_\beta^p\int_\Omega  {\abs{u}^p}{ \abs z^{ \beta}} d\xi
                 \le\int_{\Omega} \abs{\grl u}^{p} \abs z^{p+\beta} d\xi,
  \end{equation}
  where $b_{\beta}\decl\frac{\abs{m+\beta}}{p}$.
  In particular, if $N_S$ is the homogeneous norm defined in (\ref{eq:norms}),
  then for every $u\in\Cuno_0(\Omega)$, we obtain
  \begin{equation}
    \left(\frac{\abs{m-p}}{p}\right)^p\int_\Omega\frac{\abs{u}^p}{N_S^p}d\xi
    \le\left(\frac{\abs{m-p}}{p}\right)^p\int_\Omega\frac{\abs{u}^p}{ \abs z^p}d\xi
       \le \int_{\Omega} \abs{\grl u}^{p} d\xi, \label{hardycarnpsi2}
  \end{equation}
  and if $N$ is any homogeneous norm and $p\neq m$, then there exists a constant $c>0$
    such that for every $u\in\Cuno_0(\Omega)$, we have
  \begin{equation}
    c\int_\Omega\frac{\abs{u}^p}{N^p}d\xi
       \le \int_{\Omega} \abs{\grl u}^{p} d\xi, \label{hardycarnpsi3}
  \end{equation}
\item  Let $p=m>1$.
  Let $R>0$ and set $\Omega\decl\{(z,\xi_{m+1},\dots,\xi_N)\in\R^{m}\times \R^{N-m},\ \abs z <R\}$.
  If $\beta<-1$, then for every $u\in\Cuno_0(\Omega)$, we have
 \begin{equation}\label{eq:hardycarnln2}
   \tilde b_\beta^p\int_\Omega  \frac{\abs{u}^p}{\abs z^p}
       (\ln\frac R{\abs z})^\beta d\xi
     \le\int_{\Omega} \abs{\grl u}^{p} (\ln\frac R{\abs z})^{p+\beta} d\xi,
  \end{equation}
  where $\tilde b_\beta\decl\frac{\abs{\beta+1}}{p}$.
  In particular,  for every $u\in\Cuno_0(\Omega)$, we obtain
 \begin{equation}\label{dis:harcanparQ} \left(\frac{p-1}{p}\right)^p
   \int_\Omega  \frac{\abs{u}^p}{(\abs z\ln(R/\abs z))^p}d\xi
     \le\int_{\Omega} \abs{\grl u}^{p} d\xi.
  \end{equation}
\end{enumerate}
 Moreover, if $\G$ is of $H$-type, $m=l$ and
 $B^l_r\setminus\{0\}\times B^{N-l}_r\subset \Omega$
  for some $r$,
 then the constants $b_\beta^p$ and $\tilde b_\beta^p$ in
 (\ref{eq:hardycarn2}) and  in (\ref{eq:hardycarnln2})
(and hence the constants in (\ref{hardycarnpsi2}) and (\ref{dis:harcanparQ}))
are sharp.
\et

\boss From the above Theorem \ref{teo:harcarn2}, taking $m=1$, we obtain
the inequality (\ref{hardycarnpsi3}) for any $p>1$, any homogeneous norm $N$
and any function  $u\in\Cuno_0(\Omega)$ with 
$\Omega\subset(\R\setminus\{0\})\times\R^{N-1}$ and hence
also for any smooth function defined on the cone $\R_+\times\R^{N-1}$.
\eoss

\bp{} The inequalities (\ref{eq:hardycarn2})  and
(\ref{eq:hardycarnln2})
  are a direct consequence of Theorem \ref{teo:harspec}.
  The fact that $N_S\ge\abs z$ yields the inequality (\ref{hardycarnpsi2}).
  Finally the equivalence between homogeneous norms implies (\ref{hardycarnpsi3}).

  We have to prove the sharpness of the constant in the case $\G$ is of
  $H$-type with $m=l=\dim V_1$, $z=x$ and $k\decl\dim V_2$.
    We prove the optimality of the constant $b_\beta^p$ in
    (\ref{eq:hardycarn2}). The proof for $\tilde b_\beta^p$ in
    (\ref{eq:hardycarnln2}) is similar.

  We shall proceed as in the proof Theorem \ref{teo:hargr2},
  therefore it is sufficient to prove the claim for
  $\Omega=(\R^m\setminus\{0\})\times\R^{N-m}$.

  Let $c_b$ be the best constant in (\ref{eq:hardycarn2}).
  We choose $\phi\in\Cuno_0(\Omega)$ such that $\phi=uw$ with
  $u=u(\abs x)$, $w=w(t)$,
  $v\in\Cuno_0(]0,+\infty[)$ and  $w\in\Cuno_0(\R^k)$.

    Arguing as in the proof of Theorem \ref{teo:hargr2}, using the
  convexity of the function $(r^2+s^2)^{p/2}$, and the fact that
    $$\abs{\grl \phi}^2=\abs{\nabla_x \phi}^2+\frac{\abs x}{4}\abs{\nabla_t
    \phi}^2= w^2(u'(\abs x))^2+\frac{\abs x}{4}u^2\abs{\nabla_t
    w}^2,$$
  we obtain for
  $0<\lambda<1$
  \begin{equation}\label{eq:tec3}
    c_b\le (1-\lambda)^{1-p}\frac{\int_{\RR^l}\abs{\nabla_x u}^p\abs
x^{\beta+p}dx}
              {\int_{\RR^l}\abs u^p\abs x^\beta dx}
              +\lambda^{1-p}
              \frac{\int_{\RR^k}\abs{\nabla_t w}^p dt}{\int_{\RR^k}\abs w^p dt}
              \frac{ \int_{\RR^l}\abs u^p\abs x^{\gamma+\beta+p}}
               {\int_{\RR^l}\abs u^p\abs x^\beta}.
  \end{equation}
  The infimum of
  ${\int_{\RR^k}\abs{\nabla_t w}^p dt}/\int_{\RR^k}\abs w^p dt$
  vanishes.
  The infimum of
$$ \frac{ \int_{\RR^l}\abs{\nabla_xu}^p\abs x^{\gamma+\beta+p}}
               {\int_{\RR^l}\abs u^p\abs x^\beta}=
               \frac{\int_0^{+\infty} \abs{u'(s)} s^{\beta+p+l-1} ds}{\int_0^{+\infty} \abs{u(s)} s^{\beta+l-1}
               ds}$$
  is $b_\beta^p$. Indeed, it follows from Theorem \ref{teo:harcarn} with $Q=1$ and $\beta$ replaced by $\beta+l-1$.
  Letting $\lambda\rightarrow 0$ in (\ref{eq:tec3}), we conclude the proof.
\ep
\medskip

The next results deal with Hardy inequalities for functions defined on a ball
or on the complement of a ball and involving the distance from the
boundary.

If $\G$ is the Euclidean space or an $H$-type group, then the
pseudo-distance $d_2(\xi,\eta)\decl N_2(\xi^{-1}\eta) $  is actually a distance
(see \cite{cyg81}).
In a general Carnot group, there holds only the pseudo-triangular 
inequality (\ref{dis:pseudotr}).
Hence, $d_2$, in general is not a distance.
Therefore, in the general framework we shall deal with
the Carnot-Carath\'eodory distance $d_{CC}$, defined as follows.
Let $\gamma\colon[a,b]\to\RN$ be a piecewise smooth curve,
we call $\gamma$ a \emph{horizontal path} if $\dot\gamma(t)$
belongs to $V_1$ whenever it exists. 
Then for every $\xi,\eta\in\G$, we define
\be 
\label{def:dcc}
 d_{CC}(\xi,\eta)\decl\inf\left\{ \int\abs{\dot\gamma},\parbox{100mm}{ 
  $\gamma\colon[a,b]\to\RN$ horizontal path with
  $\gamma(a)=\xi,\ \gamma(b)=\eta$
}\right\}.
\ee
In the framework of Carnot group, by Chow Theorem, for every $\xi,\eta\in\G$,
it results $d_{CC}(\xi,\eta)<\infty$, and hence $d_{CC}$ is a metric on $\G$.
The distance $d_{CC}$ is left invariant with respect to the group action 
and it is homogeneous of degree 1 with respect to dilation $\delta_\lambda$, namely
$$d_{CC}(\zeta\xi,\zeta\eta)=d_{CC}(\xi,\eta),\quad
  d_{CC}(\delta_\lambda(\xi),\delta_\lambda(\eta))=\lambda d_{CC}(\xi,\eta)
  \quad\zeta,\eta,\xi\in\G, \lambda>0.$$
Hence, $d_{CC}(\cdot ,0)$ is a homogeneous norm.

\bt\label{teo:disdist}  Let $p> 1$ and let
$\Gamma_p$ be $L_p$-harmonic on $\G\setminus\{0\}$.
  Let $R>0$ and set $\Omega\decl\{ \xi\in\G, N_2(\xi)<R\}$. We have
 \begin{equation}\left(\frac{p-1}{p}\right)^p
   \int_\Omega  \frac{\abs{u}^p}{(R-N_2)^p}\abs{\grl N_2}^p d\xi
     \le\int_{\Omega} \abs{\grl u}^{p} d\xi, \qquad
     u\in D_L^{1,p}(\Omega).\label{dis:dist0}
  \end{equation}
  The constant $\left(\frac{p-1}{p}\right)^p$ is optimal.

 Moreover, we have
 \begin{equation} c^p
   \int_\Omega  \frac{\abs{u}^p}{\delta^p}\abs{\grl N_2}^p d\xi
     \le\int_{\Omega} \abs{\grl u}^{p} d\xi,\qquad u\in D_L^{1,p}(\Omega) \label{dis:dist}
  \end{equation}
  where $\delta$ is one of the following functions

   a) $\delta(\xi)\decl d_{CC}(\xi,\partial\Omega)\decl
    \inf\{d_{CC}(\xi,\eta),\ \eta\in\partial \Omega \}$,
    $c\decl \frac{p-1}{p} \frac{1}{C_2} $  and $C_2\decl\norm{\grl N_2}_{L^\infty}$;

  \noindent or

  b)  $\delta(\xi)\decl d_2(\xi,\partial\Omega)\decl
    \inf\{d_2(\xi,\eta),\ \eta\in\partial \Omega \}$ and $c\decl \frac{p-1}{p}$ 
    provided $d_2(\xi,\eta)\decl N_2(\xi^{-1}\eta)$ is a distance.
\et
\boss The constant $c^p$ in (\ref{dis:dist})
  with this generality cannot be improved. 
  Indeed, if $\G$ is the Euclidean space $\RN$ and $\grl=\nabla$, we have
  $c=\frac{p-1}{p}$, which is the best constant
  (see \cite{mat-sob97}).
\eoss
\bp{} From (\ref{id:lpn2}), and if $p\neq2$ by polarizability of $\G$, we get
  $$L_pN_2=\abs{\grl N_2}^{p-2}L_2N_2=(Q-1)\frac{\abs{\grl N_2}^{p}}{N_2}\ge0.$$
  Therefore, choosing $d(\xi)=R-N_2(\xi)$, we are in the position
  to apply Theorem \ref{th:hm} and from (\ref{dis:hargenpart})
  we get (\ref{dis:dist0}).

  Applying the scheme outlined in Remark \ref{rm:sharp},
  we obtain the optimality of the constant.

  We prove  the inequality (\ref{dis:dist}). 
  Let $\delta=d_{CC}(\cdot,\partial\Omega)$ or $\delta=d_{2}(\cdot,\partial\Omega)$,
  let $\xi\in\Omega$ be fixed and let
  $\tau\in\partial\Omega$ be a point where the minimum is attained, that is 
  $\delta(\xi)=d_{CC}(\tau,\xi)$ or $\delta(\xi)=d_{2}(\tau,\xi)$.

  First we prove the inequality (\ref{dis:dist}) in the case b).
  The inequality (\ref{dis:dist}) follows from
  (\ref{dis:dist0}) and the fact that $d_2(\cdot,\cdot)$ is a
  distance. 
  By triangular inequality, we have
  $R=d_2(\tau,0)\le d_2(\tau,\xi)+d_2(\xi,0)=\delta(\xi)+N_2(\xi)$, and hence
  we get the inequality (\ref{dis:dist}).

  We prove the case a). By the inequality
$$ \abs{N_2(\xi)-N_2(\eta)}\le \norm{\grl N_2}_{L^\infty}d_{CC}(\xi,\eta),
    \quad\mathrm{for\ every}\quad \xi,\eta\in\G, $$
  we have
  $$R-N_2(\xi)=N_2(\tau)-N_2(\xi)\le C_2 d_{CC}(\tau,\xi)=C_2\delta(\xi),$$
  which concludes the proof.
\ep

\boss We remark that if $\phi$ is a regular H-convex function and
  $\grl\phi\neq 0$ a.e., then for $p\ge2$, from (\ref{eq:lpexplic}),
  we have $-L_p(\phi)\le0$,
  thus in order to obtain Hardy inequalities involving the
  function $\phi$ we can apply the results of previous section.
  For H-convex function on Carnot groups, we refer the interested
  reader to \cite{dan-gar-nhi03,lu-man-str04,mag03}.
  For instance, in \cite{dan-gar-nhi03} the authors prove that
  in an H-type group the gauge $N$ defined in (\ref{eq:normH}) is
  H-convex, hence  $R-N(\xi)$ is H-concave and we can
  obtain again the inequalities (\ref{dis:dist0}) and (\ref{dis:dist}).
\eoss

We conclude with a Hardy inequality on an exterior domain.
\bt Let $p>Q$ and let
$\Gamma_p$ be $L_p$-harmonic on $\G\setminus\{0\}$.
  Let $R>0$ and set $\Omega\decl\{ \xi\in\G, N_2(\xi)>R\}$.
  For every $u\in \Cuno_0(\Omega)$ we have
 \begin{equation}\left(\frac{\abs{p-Q}}{p}\right)^p
   \int_\Omega  \frac{\abs{u}^p}{(N_2-R)^p}\abs{\grl N_2}^p d\xi
     \le\int_{\Omega} \abs{\grl u}^{p} d\xi.
     \label{dis:distext0}
  \end{equation}
  Moreover, for every $u\in\Cuno_0(\Omega)$, we have
 \begin{equation} c^p
   \int_\Omega  \frac{\abs{u}^p}{\delta^p}\abs{\grl N_2}^p d\xi
     \le\int_{\Omega} \abs{\grl u}^{p} d\xi,\label{dis:distext}
  \end{equation}
  where $\delta$ is one of the following functions

  a) $\delta(\xi)\decl d_{CC}(\xi,\partial\Omega)\decl
    \inf\{d_{CC}(\xi,\eta),\ \eta\in\partial \Omega \}$,
    $c\decl\frac{\abs{p-Q}}{p} \frac{1}{C_2} $  and $C_2\decl\norm{\grl N_2}_{L^\infty}$;

\noindent or

  b)
    $\delta(\xi)\decl d_2(\xi,\partial\Omega)\decl
    \inf\{d_2(\xi,\eta),\ \eta\in\partial \Omega \}$ and $c\decl\frac{\abs{p-Q}}{p}$ 
    provided $d_2(\xi,\eta)\decl N_2(\xi^{-1}\eta)$ is a distance.
\et
\bp{} Let $d$ be defined as
  $d(\xi)\decl N_2(\xi)^{\frac{p-Q}{p-1}}-R^{\frac{p-Q}{p-1}}$
  ($\xi\in\Omega$). It is clear that $d$ is positive and
  $L_p d=0$. Applying Theorem \ref{th:hm} we derive
  \begin{equation}\left(\frac{\abs{p-Q}}{p}\right)^p
   \int_\Omega {\abs{u}^p}
   \left( \frac{N_2^{\frac{1-Q}{p-1}}}{N_2(\xi)^{\frac{p-Q}{p-1}}-R^{\frac{p-Q}{p-1}}}\right)^p \abs{\grl N_2}^p d\xi
     \le\int_{\Omega} \abs{\grl u}^{p} d\xi.
     \label{dis:distext00}
  \end{equation}
  It is easy to check that for $\xi\in\Omega$, it results
  $$N_2(\xi)^{\frac{p-Q}{p-1}}-R^{\frac{p-Q}{p-1}}\le
    {N_2^{\frac{1-Q}{p-1}}}(N_2-R), $$
  which with (\ref{dis:distext00}) implies (\ref{dis:distext0}).

  Arguing as in the proof of Theorem \ref{teo:disdist},
  we obtain the missing inequality (\ref{dis:distext}).
\ep

\boss The constant $c^p$ in
 (\ref{dis:distext0}) and (\ref{dis:distext}) cannot be improved
 in this generality. Indeed if $\grl$ is the usual
 gradient $\nabla$, then this constant is sharp
 (see \cite{mat-sob98}).
\eoss


\def\cprime{$'$} \def\cprime{$'$} \def\cprime{$'$} \def\cprime{$'$}
  \def\cprime{$'$} \def\cprime{$'$}

\end{document}